\documentclass[12pt]{article}
\usepackage{amssymb, amsmath, amsthm, amscd}
\usepackage[utf8]{inputenc}
\usepackage[all,cmtip]{xy}
\usepackage{enumitem}
\usepackage{setspace}
\usepackage{mathdots}

\usepackage[mathscr]{euscript}
\usepackage[colorlinks=true]{hyperref}
\hypersetup{colorlinks   = true}
\hypersetup{linkcolor=blue}
\usepackage{tikz}

\begin{document}

\mathchardef\mhyphen="2D
\newtheorem{The}{Theorem}[section]
\newtheorem{Lem}[The]{Lemma}
\newtheorem{Prop}[The]{Proposition}
\newtheorem{Cor}[The]{Corollary}
\newtheorem{Rem}[The]{Remark}
\newtheorem{Obs}[The]{Observation}
\newtheorem{SConj}[The]{Standard Conjecture}
\newtheorem{Titre}[The]{\!\!\!\! }
\newtheorem{Conj}[The]{Conjecture}
\newtheorem{Question}[The]{Question}
\newtheorem{Prob}[The]{Problem}
\newtheorem{Def}[The]{Definition}
\newtheorem{Not}[The]{Notation}
\newtheorem{Claim}[The]{Claim}
\newtheorem{Conc}[The]{Conclusion}
\newtheorem{Ex}[The]{Example}
\newtheorem{Fact}[The]{Fact}
\newtheorem{Formula}[The]{Formula}
\newtheorem{Formulae}[The]{Formulae}
\newtheorem{The-Def}[The]{Theorem and Definition}
\newtheorem{Prop-Def}[The]{Proposition and Definition}
\newtheorem{Lem-Def}[The]{Lemma and Definition}
\newtheorem{Cor-Def}[The]{Corollary and Definition}
\newtheorem{Conc-Def}[The]{Conclusion and Definition}
\newtheorem{Terminology}[The]{Note on terminology}
\newcommand{\C}{\mathbb{C}}
\newcommand{\R}{\mathbb{R}}
\newcommand{\N}{\mathbb{N}}
\newcommand{\Z}{\mathbb{Z}}
\newcommand{\Q}{\mathbb{Q}}
\newcommand{\Proj}{\mathbb{P}}
\newcommand{\Rc}{\mathcal{R}}
\newcommand{\Oc}{\mathcal{O}}
\newcommand{\Vc}{\mathcal{V}}
\newcommand{\Id}{\operatorname{Id}}
\newcommand{\pr}{\operatorname{pr}}
\newcommand{\rk}{\operatorname{rk}}
\newcommand{\del}{\partial}
\newcommand{\delbar}{\bar{\partial}}
\newcommand{\Cdot}{{\raisebox{-0.7ex}[0pt][0pt]{\scalebox{2.0}{$\cdot$}}}}
\newcommand\nilm{\Gamma\backslash G}
\newcommand\frg{{\mathfrak g}}
\newcommand{\fg}{\mathfrak g}
\newcommand{\Oh}{\mathcal{O}}
\newcommand{\Kur}{\operatorname{Kur}}
\newcommand\gc{\frg_\mathbb{C}}
\newcommand\hisashi[1]{{\textcolor{red}{#1}}}
\newcommand\dan[1]{{\textcolor{blue}{#1}}}
\newcommand\luis[1]{{\textcolor{orange}{#1}}}

\begin{center}

{\Large\bf Properties of Holomorphic $p$-Contact Manifolds}

\end{center}

\begin{center}

{\large Hisashi Kasuya, Dan Popovici and Luis Ugarte}

\end{center}

\vspace{1ex}

\noindent{\small{\bf Abstract.} We continue the study of compact holomorphic $p$-contact manifolds $X$ that we introduced recently by expanding the discussion to include non-K\"ahler hyperbolicity issues and a differential calculus based on what we call the Lie derivative with respect to a $(0,\,q)$-form with values in the holomorphic tangent bundle of $X$. We also propose the notion of $p$-contact deformations for which we prove a Bogomolov-Tian-Todorov-type unobstructedness theorem to order two. This kind of small deformations of the complex structure is related to the essential horizontal deformations that we introduced in our previous work and forms part of a wider on-going project aimed at developing a non-K\"ahler mirror symmetry theory that was first tested on the Iwasawa manifold and subsequently on Calabi-Yau page-$1$-$\partial\bar\partial$-manifolds.}

\vspace{1ex}


\vspace{1ex}


\vspace{1ex}

\section{Introduction}\label{section:introduction}

In the first part of this work, we have recently introduced the class of compact complex manifolds described in the following

\begin{Def}([KPU25])\label{Def:hol_p-contact_structure} Let $X$ be a compact complex manifold with $\mbox{dim}_\C X =n = 2p+1$.

\vspace{1ex}

\noindent $(1)$\, A {\bf holomorphic $p$-contact structure} on $X$ is a smooth $(p,\,0)$-form $\Gamma\in C^\infty_{p,\,0}(X,\,\C)$ such that \begin{eqnarray*}(a)\hspace{2ex} \bar\partial\Gamma=0 \hspace{6ex} \mbox{and} \hspace{6ex} (b)\hspace{2ex} \Gamma\wedge\partial\Gamma\neq 0 \hspace{2ex} \mbox{at every point of} \hspace{1ex} X.\end{eqnarray*}

\vspace{1ex}

\noindent $(2)$\, We say that $X$ is a {\bf holomorphic $p$-contact manifold} if there exists a  holomorphic $p$-contact structure $\Gamma$ on $X$.

\end{Def}

\vspace{1ex}

In this second part, we continue the study of this notion by pointing out new examples, developing a notion of Lie derivative and investigating the hyperbolicity properties and a different type of small deformations of the complex structures of holomorphic $p$-contact manifolds.

\vspace{2ex}

 (I)\, We will sometimes contrast our notion of {\it holomorphic $p$-contact structure} with a kind of opposite, reminiscent of [Dem02], introduced in

\begin{Def}\label{Def:hol_p-no-contact_structure} Let $X$ be a compact complex manifold with $\mbox{dim}_\C X =n = 2p+1$.

 \vspace{1ex}

\noindent $(1)$\, A {\bf holomorphic $p$-no-contact structure} on $X$ is a form $\Gamma\in C^\infty_{p,\,0}(X,\,\C)$ such that \begin{eqnarray*}(a)\hspace{2ex} \bar\partial\Gamma=0 \hspace{6ex} \mbox{and} \hspace{6ex} (b)\hspace{2ex} \partial\Gamma = \Gamma\wedge\zeta \hspace{2ex} \mbox{for some form} \hspace{1ex} \zeta\in C^\infty_{1,\,0}(X,\,\C).\end{eqnarray*}

\vspace{1ex}

\noindent $(2)$\, We say that $X$ is a {\bf holomorphic $p$-no-contact manifold} if there exists a  holomorphic $p$-no-contact structure $\Gamma$ on $X$.

\end{Def}  

Note that, if $p$ is {\it odd}, the condition $\partial\Gamma = \Gamma\wedge\zeta$ implies $\Gamma\wedge\partial\Gamma = (\Gamma\wedge\Gamma)\wedge\zeta = 0$. In this sense, a holomorphic $p$-no-contact structure is the opposite of a holomorphic $p$-contact one.

\vspace{1ex}

(II)\, A key technical role in this paper is played by a differential calculus that we introduce in $\S$\ref{section:Lie}. (A similar approach was used by [Cle05] in a different context.) It involves what we call the {\it Lie derivative} with respect to a vector-valued form $\theta\in C^{\infty}_{0,\,q}(X,\,T^{1,\,0}X)$ on a complex manifold $X$ (see Definition \ref{Def:Lie-derivatives-mixed_def}), the case $q=1$ being the most of interest to us. This Lie derivative is modelled on the classical Lie derivative with respect to a vector field and is motivated by the need for a general formula computing expressions of the shape $\partial(\theta\lrcorner u)$, where $u$ is a scalar-valued differential form of any (bi-)degree and the contraction operation $\theta\lrcorner\cdot$ combines the contraction of $u$ by the vector-field part of $\theta$ with the exterior product by the $(0,\,q)$-form part of $\theta$. Such a formula is indispensable in handling the ubiquitous quantities of the form $\theta\lrcorner u$ arising on Calabi-Yau manifolds of which our holomorphic $p$-contact manifolds $(X,\,\Gamma)$ are a subclass. 

Specifically, with every $\theta\in C^{\infty}_{0,\,q}(X,\,T^{1,\,0}X)$, we associate the linear differential operator $L_\theta:C^\infty_{r,\,s}(X,\,\C)\longrightarrow C^\infty_{r,\,s+q}(X,\,\C)$ of order one and bidegree $(0,\,q)$ defined so that the following analogue of the Leibniz formula holds in this case: \begin{equation*}\partial(\theta\lrcorner u) = L_\theta u + (-1)^{q-1}\,\theta\lrcorner\partial u.\end{equation*} The above sign is due to the total degree of any such $\theta$ being $q-1$, with the $T^{1,\,0}X$-part counting for $-1$. The case $q=0$ was introduced and used in [PU23], where the Lie derivative operator with respect to a $(1,\,0)$-vector field $\xi\in C^{\infty}(X,\,T^{1,\,0}X)$ was denoted by $L^{1,\,0}_{\xi} :C^\infty_{r,\,s}(X,\,\C)\longrightarrow C^\infty_{r,\,s}(X,\,\C)$.

In Proposition \ref{Prop:Lie-derivatives_local} and Lemma \ref{Lem:Lie-derivatives-mixed_prop} we give the main properties, including various commutation identities, of the Lie derivatives $L_\theta$.

Applications of this Lie derivative calculus include the emergence in $\S$\ref{subsection:applications_Lie} of the notions of {\it constantly horizontal} and {\it constantly vertical} $T^{1,\,0}X$-valued $(0,\,q)$-forms on holomorphic $p$-contact manifolds. These sharpen the notions of {\it horizontal} (i.e. ${\cal F}_\Gamma$-valued), respectively {\it vertical} (i.e. ${\cal G}_\Gamma$-valued), $(0,\,q)$-forms and play a significant role in the theory of small deformations of compact holomorphic $p$-contact manifolds $(X,\,\Gamma)$ that we initiate in $\S$\ref{section:contact-def}. An equally important role is played by the various generalisations of the classical Tian-Todorov formula that we obtain in $\S$\ref{subsection:generalised_Tian-Todorov} as further generalisations of our Lie derivative calculus.

\vspace{1ex}

(III)\, Besides the various examples of compact {\it holomorphic $p$-contact manifolds} that we give in $\S$\ref{section:examples_dimension_3} and the development of a Lie derivative calculus that occupies $\S$\ref{section:Lie}, we obtain general hyperbolicity results for our manifolds. Specifically, in $\S$\ref{section:contact-hyperbolicity} we discuss the implications for the setting of holomorphic $p$-contact manifolds of the non-K\"ahler hyperbolicity theory initiated in [MP22a], [MP22b] and [KP23] and continued in [Mar23a], [Mar23b] and [Ma24].

Recall the two main peculiar features of this new point of view on hyperbolicity:


(a)\, the theory accommodates possibly {\bf non-K\"ahler} compact complex manifolds, while the classical notions of hyperbolicity in the senses of Kobayashi, Brody, Gromov, etc, apply only to projective manifolds;


(b)\, the theory rules out the existence of non-trivial entire holomorphic maps from possibly {\bf higher-dimensional} spaces $\C^p$ (i.e. possibly with $p\geq 2$) into a given compact complex manifold $X$, while the classical Brody hyperbolicity is concerned solely with maps from $\C$ to $X$. The main extra difficulty in this new setting is that no {\it Brody Renormalisation Lemma} exists, a fact that necessitated the introduction of a notion of slow growth (called {\it subexponential growth} in [MP22a] and the subsequent works) for the entire holomorphic maps $f:\C^p\longrightarrow X$ whose existence is ruled out for $X$ to qualify as hyperbolic.

\vspace{1ex}

In the particular context of this paper, we introduce two notions of hyperbolicity (cf. Definitions \ref{Def:p-contact-metric-hyperbolic} and \ref{Def:p-contact-hyperbolic}) corresponding to the {\it metric} (generalising Gromov's K\"ahler hyperbolicity), respectively {\it entire-map-based} (generalising the Brody hyperbolicity) side of the theory. That the former hyperbolicity notion implies the latter is ensured by Theorem \ref{The:partial-p-hyperbolicity_implication}, while a general result on the existence of Ahlfors currents from [KP23] is applied in Theorem \ref{The:partial-p+1-hyperbolicity_p-contact} to prove yet another hyperbolicity property that compact holomorphic $p$-contact manifolds $(X,\,\Gamma)$ have under certain extra assumptions. This result rules out the existence of holomorphic maps $f:\C^{p+1}\longrightarrow X$, while Theorem \ref{The:partial-p-hyperbolicity_implication} rules out the existence of holomorphic maps $f:\C^p\longrightarrow X$, each with a different type of slow growth, but both with the property $f^\star\Gamma = 0$. This last condition is a way of expressing a type of {\it partial} (i.e. holding only in some, not necessarily all, directions) hyperbolicity peculiar to the $p$-contact setting and different from the one studied in [KP23].

\vspace{1ex}

(IV)\, Another type of general conceptual results that we obtain in this paper pertains to small deformations of holomorphic $p$-contact manifolds $(X,\,\Gamma)$ and occupies $\S$\ref{section:contact-def}. In Definition \ref{Def:contact_01_T10}, we propose the notion of {\bf $p$-contact deformations} of $(X,\,\Gamma)$ which we then prove in Theorem \ref{The:partial-unobstructedness} to be what we call {\bf unobstructed to order two}, a partial unobstructedness property.

This notion bears a similarity to the notion of {\it essential horizontal deformations} introduced in [KPU25] as Definition 6.1. Each has its own role to play, so we present them separately. In particular, we believe that they will both contribute to the development of a non-K\"ahler mirror symmetry theory that started in [Pop18] with the Iwasawa manifold and continued in [PSU20a], [PSU20b] and [PSU20c].

\section{Preliminaries}\label{section:preliminaries} We will make two kinds of observations about these classes of manifolds.

\subsection{Incompatibility with certain currents}\label{subsection:incompatibility} The next observation implies that the existence of a strongly (strictly) positive $\partial\bar\partial$-closed current $T$ of bidegree $(p,\,p)$ on $X$ rules out the existence of a holomorphic $p$-contact structure.

The notion of a current $T$ of bidegree $(p,\,p)$ being {\it strongly (strictly) positive} on a compact $n$-dimensional complex manifold $X$ was introduced in $(1)$ of Definition $3.25$ of the appendix to [Pop25]. It means that, having fixed any Hermitian metric $\omega$ on $X$, there exists a constant $\varepsilon>0$ such that the $(p,\,p)$-current $T - \varepsilon\omega^p$ is {\it strongly semi-positive} on $X$ (in the usual sense that for any weakly semi-positive $C^\infty$ $(n-p,\,n-p)$-form $\Omega$ on $X$, we have $\langle T - \varepsilon\omega^p,\,\Omega\rangle \geq 0$.) We write in this case $T>0$ (strongly) or $T\geq\varepsilon\omega^p$ (strongly) on $X$. Thanks to the compactness of $X$, if such a constant $\varepsilon>0$ exists for some Hermitian metric $\omega$ on $X$, it also exists for any Hermitian metric.

Similarly, a $C^\infty$ $(p,\,p)$-form $\Xi$ on $X$ is said to be {\it strongly (strictly) positive} (cf. $(2)$ of Definition $3.11$ of [Pop25]) if $\Xi - \varepsilon\omega^p$ is a {\it strongly semi-positive} form on $X$ for some constant $\varepsilon>0$ and some Hermitian metric $\omega$. Again, we write in this case $\Xi>0$ (strongly) or $\Xi\geq\varepsilon\omega^p$ (strongly) on $X$.

The notion of {\it strongly positive non-degenerate} $(p,\,p)$-form at a pregiven point was introduced in [HK74, Corollary 1.10] by requiring that it belong to the interior of the classical {\it cone of strongly semi-positive (p, p)-forms}, denoted by $SP^p$ in [HK74] and by ${\cal C}$ in [Pop25]. This notion turns out to be equivalent to the notion of {\it strongly (strictly) positive} $(p,\,p)$-form at a pregiven point (cf. Remark 3.8 in [Pop25]). However, it is the version for currents and forms on a compact manifold, recalled above, that will be used in this section.   

Recall that any $C^\infty$ $(p,\,p)$-form defines a current and that the semi-positivity (in either the strong or the weak sense) of such a form in the sense of forms is equivalent to the analogous notion in the sense of currents -- see, for example, Observation $3.24$ in [Pop25]. Hence, \begin{eqnarray*}\Xi>0 \hspace{2ex} \mbox{(strongly) \hspace{2ex} as a form} \iff  \Xi>0 \hspace{2ex} \mbox{(strongly) \hspace{2ex} as a current}.\end{eqnarray*} 

For an overview of the classical notions of semi-positivity introduced by Lelong in the 1950s, see [Dem97, III, $\S1$] or the appendix to [Pop25] (where semi-positivity is reviewed before (strict) positivity is introduced and investigated). In particular, recall that a $(p,\,p)$-form $\Xi$ is said to be {\it strongly semi-positive} at a point $x\in X$ if it is a {\it finite} linear combination with {\it non-negative} coefficients $\lambda_s \geq 0$ of decomposable $(p,\,p)$-forms: \begin{eqnarray*}\Xi = \sum\limits_s \lambda_s\, i\alpha_{s,\,1}\wedge\overline\alpha_{s,\,1}\wedge\dots\wedge i\alpha_{s,\,p}\wedge\overline\alpha_{s,\,p},\end{eqnarray*} where $\alpha_{s,\,1}, \dots , \alpha_{s,\,p}$ are $(1,\,0)$-forms at $x$. Further recall that a $(p,\,p)$-form $\Xi$ is said to be {\it weakly semi-positive} at a point $x\in X$ if it evaluates {\it non-negatively} against all decomposable forms of the complementary bidegree $(n-p,\,n-p)$. This means that \begin{eqnarray*}\Xi\wedge i\alpha_1\wedge\overline\alpha_1\wedge\dots\wedge i\alpha_{n-p}\wedge\overline\alpha_{n-p}\end{eqnarray*} is a {\it non-negative} $(n,\,n)$-form for all families $\{\alpha_1,\dots , \alpha_{n-p}\}$ of $(1,\,0)$-forms at $x$. 

\begin{Obs}\label{Obs:incompatibility} Let $X$ be a compact complex manifold with $\mbox{dim}_\C X = n = 2p+1$. A {\bf holomorphic $p$-contact structure} $\Gamma\in C^\infty_{p,\,0}(X,\,\C)$ and a {\bf strongly (strictly) positive $\partial\bar\partial$-closed current} $T$ of bidegree $(p,\,p)$ cannot simultaneously exist on $X$.

  In particular, no $3$-dimensional compact complex manifold that carries an {\bf SKT metric} can carry a holomorphic contact structure. Similarly, in arbitrary dimension, no compact complex manifold that carries a {\bf strongly (strictly) positive $\partial\bar\partial$-closed} $C^\infty$ $(p,\,p)$-form can carry a holomorphic $p$-contact structure.

\end{Obs}  

\noindent {\it Proof.} Suppose there exists a holomorphic $p$-contact structure $\Gamma\in C^\infty_{p,\,0}(X,\,\C)$. By Remark 1.10 in [Dem97, III], the $(p+1,\,p+1)$-form $i^{(p+1)^2}\partial\Gamma\wedge\bar\partial\overline\Gamma$ is weakly semi-positive. It is also $\partial\bar\partial$-exact since \begin{eqnarray*}i^{(p+1)^2}\partial\Gamma\wedge\bar\partial\overline\Gamma = \partial\bar\partial(-i^{(p+1)^2}\Gamma\wedge\overline\Gamma).\end{eqnarray*}

If a strongly (strictly) positive $\partial\bar\partial$-closed current $T$ of bidegree $(p,\,p)$ existed on $X$, then $i^{(p+1)^2}\partial\Gamma\wedge\bar\partial\overline\Gamma\wedge T$ would be a non-negative $(n,\,n)$-current on $X$, so we would get the inequality below: \begin{eqnarray*}0\leq\int_X i^{(p+1)^2}\partial\Gamma\wedge\bar\partial\overline\Gamma\wedge T = -i^{(p+1)^2} \int_X \Gamma\wedge\overline\Gamma\wedge\partial\bar\partial T = 0,\end{eqnarray*} where the first equality follows from Stokes's theorem and the second equality follows from $\partial\bar\partial T = 0$. This would imply that $i^{(p+1)^2}\partial\Gamma\wedge\bar\partial\overline\Gamma = 0$, hence also that \begin{eqnarray*}dV_\Gamma := i^{n^2}\,(\Gamma\wedge\partial\Gamma)\wedge(\overline\Gamma\wedge\bar\partial\overline\Gamma) =0.\end{eqnarray*}

However, this is a contradiction since the $C^\infty$ volume form $dV_\Gamma$ is strictly positive at every point of $X$ thanks to the hypothesis (implicit in $\Gamma$ being a holomorphic $p$-contact structure) that $\Gamma\wedge\partial\Gamma$ is non-vanishing everywhere on $X$.  \hfill $\Box$

\subsection{Integrability}\label{subsection:integrability} The next observation, in the vein of [Dem02], is that any {\it holomorphic $p$-no-contact structure} induces a {\it (possibly singular) foliation}, namely an {\it integrable} coherent (possibly non-locally free) subsheaf ${\cal F}_\Gamma$ of the holomorphic tangent sheaf ${\cal O}(T^{1,\,0}X)$.

\begin{Obs}\label{Obs:no-contact_foliation} Suppose there exists a {\bf holomorphic $p$-no-contact structure} $\Gamma\in C^\infty_{p,\,0}(X,\,\C)$ on a compact complex $(2p+1)$-dimensional manifold $X$.

  Then, the coherent subsheaf ${\cal F}_\Gamma$ of ${\cal O}(T^{1,\,0}X)$ of germs of $(1,\,0)$-vector fields $\xi$ such that $\xi\lrcorner\Gamma=0$ is {\bf integrable} in the sense that $[{\cal F}_\Gamma,\,{\cal F}_\Gamma]\subset{\cal F}_\Gamma$, where $[\,\cdot\,,\,\cdot\,]$ is the Lie bracket of $T^{1,\,0}X$.

\end{Obs}

\noindent {\it Proof.} Since $\bar\partial\Gamma=0$, we have $\partial\Gamma = d\Gamma$. Meanwhile, the well-known Cartan formula (whose use is standard in this type of contexts, see e.g. [Dem02]) yields the latter equality below: \begin{eqnarray}\label{eqn:Cartan-formula}\nonumber (\Gamma\wedge\zeta)(\xi_0,\dots ,\xi_p) = (d\Gamma)(\xi_0,\dots ,\xi_p) & = & \sum\limits_{j=0}^p(-1)^j\,\xi_j\cdot\Gamma(\xi_0,\dots ,\widehat\xi_j,\dots , \xi_p) \\
  & +  & \sum\limits_{0\leq j<k\leq p}(-1)^{j+k}\,\Gamma([\xi_j,\,\xi_k],\,\xi_0,\dots,\widehat\xi_j,\dots ,\widehat\xi_k,\dots , \xi_p)\end{eqnarray} for all $(1,\,0)$-vector fields $\xi_0,\dots , \xi_p$.

If two of the vector fields $\xi_0,\dots ,\xi_p$, say $\xi_l$ and $\xi_s$ for some $l<s$, are (local) sections of ${\cal F}_\Gamma$, we have:

\vspace{1ex}

(a)\, $(\Gamma\wedge\zeta)(\xi_0,\dots ,\xi_p) = 0$ since $\zeta$ is a $1$-form, so it can absorb at most one of $\xi_l$ and $\xi_s$, say $\xi_l$, hence at least one of them, $\xi_s$ in this case, is an argument of $\Gamma$ and we can use the property $\xi_s\lrcorner\Gamma = 0$;

\vspace{1ex}

(b)\, all the terms in the first sum on the r.h.s. of (\ref{eqn:Cartan-formula}) vanish;

 \vspace{1ex}

(c)\, all the terms in the second sum on the r.h.s. of (\ref{eqn:Cartan-formula}) vanish, except possibly the term $$(-1)^{l+s}\,\Gamma([\xi_l,\,\xi_s],\,\xi_0,\dots,\widehat\xi_l,\dots ,\widehat\xi_s,\dots , \xi_p).$$

 \vspace{1ex}

 But then this last term must vanish as well for all choices of $(1,\,0)$-vector fields $\xi_0,\dots,\widehat\xi_l,\dots ,\widehat\xi_s,\dots , \xi_p$. This means that $[\xi_l,\,\xi_s]\lrcorner\Gamma = 0$, which amounts to $[\xi_l,\,\xi_s]$ being a (local) section of ${\cal F}_\Gamma$, for all pairs $\xi_l,\,\xi_s$ of (local) sections of ${\cal F}_\Gamma$.

 This proves the integrability of ${\cal F}_\Gamma$.  \hfill $\Box$

 \section{Examples in dimension $3$}\label{section:examples_dimension_3} We will present five 
classes of examples that have different properties w.r.t. our notions. Moreover, we discuss the uniqueness of these examples in complex dimension 3 in the context of Lie groups endowed with left invariant complex structures and their quotients by lattices.

 \vspace{1ex}

$(1)$\, The {\bf Iwasawa manifold} $X=G/\Lambda$ is the quotient of the nilpotent complex Lie group (called the Heisenberg group) $G=(\C^3,\,\star)$ whose group operation is defined by $$(\zeta_1,\,\zeta_2,\,\zeta_3)\star(z_1,\,z_2,\,z_3)=(\zeta_1+z_1,\,\zeta_2+z_2,\,\zeta_3+z_3+\zeta_1\,z_2),$$ by the lattice $\Lambda\subset G$ consisting of the elements $ (z_1,\,z_2,\,z_3)\in G$ with $z_1,\,z_2,\,z_3\in\Z[i]$. (See e.g. [Nak75].)

 The holomorphic $(1,\,0)$-forms $dz_1,\,dz_2,\, dz_3-z_1\,dz_2$ on $\C^3$ induce {\it holomorphic} $(1,\,0)$-forms $\alpha,\,\beta,\,\gamma$ on $X$ that satisfy the structure equations: $$\partial\alpha = \partial\beta = 0  \hspace{5ex}\mbox{and}\hspace{5ex} \partial\gamma = -\alpha\wedge\beta\neq 0 \hspace{1ex}\mbox{at every point of}\hspace{1ex} X.$$

 Thus, $\Gamma:=\gamma$ defines a {\bf holomorphic contact structure} ($=$ a {\bf holomorphic $1$-contact structure)} since $$\gamma\wedge\partial\gamma =  -\alpha\wedge\beta\wedge\gamma \neq 0 \hspace{2ex}\mbox{at every point of}\hspace{1ex} X.$$ This fact is actually standard.

 On the other hand, each of the forms $\alpha$ and $\beta$ defines a {\it holomorphic $1$-no-contact structure}, but these structures are trivial since $\partial\alpha = \partial\beta = 0$. 

 \vspace{1ex}

 $(2)$\, The {\bf Nakamura manifolds} $X=G/\Lambda$ are quotients of a solvable, non-nilpotent complex Lie group $G=(\C^3,\,\star)$ whose group operation is defined by $$(\zeta_1,\,\zeta_2,\,\zeta_3)\star(z_1,\,z_2,\,z_3)=(\zeta_1+z_1,\,\zeta_2+e^{-\zeta_1}z_2,\,\zeta_3+e^{\zeta_1}z_3),$$ by a lattice $\Lambda\subset G$. (See e.g. [Nak75].) The cohomology of any such $X$ is determined by three {\it holomorphic} $(1,\,0)$-forms $\varphi_1, \varphi_2, \varphi_3$ that satisfy the structure equations: $$\partial\varphi_1 = 0, \hspace{5ex} \partial\varphi_2 = \varphi_1\wedge\varphi_2   \hspace{5ex}\mbox{and}\hspace{5ex} \partial\varphi_3 = -\varphi_1\wedge\varphi_3.$$ Thus, each of the forms $\varphi_1, \varphi_2, \varphi_3$ defines a {\bf holomorphic $1$-no-contact structure}, those defined by $\varphi_2$ and $\varphi_3$ being non-trivial (since $\partial\varphi_j\neq 0$ for $j=2,3$.)

 On the other hand, each of the forms $\Gamma_1:= \varphi_2 + \varphi_3$ and $\Gamma_2:= \varphi_2 - \varphi_3$ defines a {\bf holomorphic contact structure} ($=$ a {\bf holomorphic $1$-contact structure)} on $X$. Indeed, we have: \begin{eqnarray*}\Gamma_1\wedge\partial\Gamma_1 = 2\,\varphi_1\wedge\varphi_2\wedge\varphi_3 \hspace{3ex} \mbox{and} \hspace{3ex} \Gamma_2\wedge\partial\Gamma_2 = -2\,\varphi_1\wedge\varphi_2\wedge\varphi_3,\end{eqnarray*} while $\varphi_1\wedge\varphi_2\wedge\varphi_3$ is a non-vanishing holomorphic $(3,\,0)$-form on $X$.

\vspace{1ex}

$(3)$\, The compact complex $3$-dimensional manifolds $X=SL(2,\,\C)/\Lambda$ that are quotients of the simple complex Lie group $SL(2,\,\C)$ by any lattice have their cohomology determined by three {\it holomorphic} $(1,\,0)$-forms $\alpha,\,\beta,\,\gamma$ that satisfy the structure equations: $$\partial\alpha = \beta\wedge\gamma, \hspace{5ex} \partial\beta = \gamma\wedge\alpha, \hspace{5ex} \partial\gamma = \alpha\wedge\beta.$$ Thus, each of the forms $\alpha,\,\beta,\,\gamma$ defines a {\bf holomorphic contact structure} ($=$ a {\bf holomorphic $1$-contact structure)} on $X$.

\vspace{1ex}

$(4)$\, {\bf A non complex parallelisable nilmanifold.} Note that the examples in (1)--(3) are quotients of complex Lie groups. Here we consider instead a (real) nilpotent Lie group $G$ endowed with a left invariant complex structure defined by complex coordinates $(z_1,z_2,z_3)$ with respect to which $G=(\C^3,\,\star)$ is given by the following (non holomorphic) group operation:
$$(\zeta_1,\,\zeta_2,\,\zeta_3)\star(z_1,\,z_2,\,z_3) = \bigg(\zeta_1+z_1,\,\zeta_2+z_2 + \overline{\zeta}_1\,z_1,\,\zeta_3+z_3-\zeta_1\,z_2 - \overline{\zeta}_1\,z_1(\zeta_1+\frac{z_1}{2})\bigg).$$
We can consider the lattice $\Lambda\subset G$ consisting of the elements $ (z_1,\,z_2,\,z_3)\in G$ with $z_1\in 2\Z[i]$ and $z_2,\,z_3\in\Z[i]$ to get a compact complex manifold $X=G/\Lambda$. 
The basis of left invariant $(1,\,0)$-forms 
$$dz_1,\ dz_2-\overline{z}_1\,dz_1,\ dz_3+z_1\,dz_2$$ 
on $\C^3$ induces a basis of $(1,\,0)$-forms $\varphi_1,\,\varphi_2,\,\varphi_3$ on $X$ that satisfy the structure equations:
 $$d \varphi_1=0,  \hspace{5ex} d \varphi_2 = \varphi_1\wedge\overline{\varphi}_1, \hspace{5ex} d \varphi_3 = \varphi_1\wedge\varphi_2.$$
Alternatively, these equations correspond to a left invariant complex structure on the nilpotent Lie group with Lie algebra ${\mathfrak h}_{15}$, whereas the example in (1) is related to the Lie algebra ${\mathfrak h}_{5}$, which is the real Lie algebra underlying the Iwasawa manifold (see [Uga07] for more details).

The compact complex manifold 
$X$ has a {\bf holomorphic contact structure} ($=$ a {\bf holomorphic $1$-contact structure)} defined by the form $\varphi_3$, since $\delbar\varphi_3=0$ and $\varphi_3\wedge\partial\varphi_3 = \varphi_1\wedge\varphi_2\wedge\varphi_3 \neq 0$.

On the other hand, the closed form $\varphi_1$ defines trivially a {\it holomorphic $1$-no-contact structure on~$X$}.

\vspace{1ex}

$(5)$\, {\bf A non-complex parallelisable solvmanifold.} 
The complex parallelisable example given in (2) corresponds to a left invariant complex structure living on the solvable Lie algebra underlying the Nakamura manifold. This Lie algebra, denoted by ${\mathfrak g}_8$ in [FOU15], has many other (non-bi-invariant) complex structures. 
Here we consider the one defined by the basis of $(1,\,0)$-forms 
$\varphi_1,\,\varphi_2,\,\varphi_3$ satisfying:
 $$d \varphi_1=0,  
 \hspace{5ex} d \varphi_2 = \varphi_1\wedge\varphi_2 +\varphi_1\wedge\overline{\varphi}_1, 
 \hspace{5ex} d \varphi_3 = - \varphi_1\wedge\varphi_3 + \varphi_1\wedge\overline{\varphi}_1.$$
 Now, if $X$ is a  compact complex solvmanifold obtained from these equations, then
$X$ has a {\bf holomorphic contact structure} ($=$ a {\bf holomorphic $1$-contact structure)} defined by the form $\Gamma=\varphi_2-\varphi_3$, since $\delbar\Gamma=0$ and $\partial\Gamma =\varphi_1\wedge(\varphi_2 +\varphi_3)$, so $\Gamma\wedge\partial\Gamma = -2\varphi_1\wedge\varphi_2\wedge\varphi_3 \neq 0$.
On the other hand, the closed form $\varphi_1$ defines trivially a {\it holomorphic $1$-no-contact structure on~$X$}.

\begin{Obs}\label{Obs:non-Iwasawa-non-page-1-nor-balanced} 
The  compact complex manifold $X$ given in $(4)$ is not a page-$1$-$\partial\bar\partial$-manifold. 
Indeed, the $(0,\,1)$-form $\overline{\varphi}_2$ satisfies the following conditions:
$$
\delbar \overline{\varphi}_2 =0,\  \ \ \partial\overline{\varphi}_2 = -\varphi_1\wedge\overline{\varphi}_1 = \delbar (-\varphi_2),
$$
so $\varphi_1\wedge\overline{\varphi}_1 \in \partial\Big({\cal Z}_2^{0,\,1}(X)\Big)$. Suppose that $\partial \overline{\varphi}_2$ belongs to $\mbox{Im}\,(\partial\bar\partial)$, i.e. $\varphi_1\wedge\overline{\varphi}_1=\partial\delbar f$ for some function $f$ on $X$. Then, $\partial\delbar f \wedge\varphi_2\wedge\overline{\varphi}_2\wedge\varphi_3\wedge\overline{\varphi}_3$ is a nowhere vanishing $(3,3)$-form on $X$, and we arrive at the contradiction:
$$
0\not=\int_X \partial\delbar f \wedge\varphi_2\wedge\overline{\varphi}_2\wedge\varphi_3\wedge\overline{\varphi}_3 = \int_X  f \, \partial\delbar(\varphi_2\wedge\overline{\varphi}_2\wedge\varphi_3\wedge\overline{\varphi}_3) =0,
$$
where the latter equality is a consequence of the fact that $\partial \delbar(\varphi_2\wedge\overline{\varphi}_2\wedge\varphi_3\wedge\overline{\varphi}_3)=\partial (\varphi_1\wedge\overline{\varphi}_1\wedge\overline{\varphi}_2\wedge\varphi_3\wedge\overline{\varphi}_3)=0$.
Thus, $\partial\Big({\cal Z}_2^{0,\,1}(X)\Big) \not\subset \mbox{Im}\,(\partial\bar\partial)$, so \eqref{eqn:ddbar_im_Z-2}~is not satisfied.

On the other hand, the  compact complex manifold $X$ given in $(4)$ admits no balanced metric. Recall that a Hermitian metric $\omega$ on an $n$-dimensional  complex manifold is said to be balanced if $d\omega^{n-1}=0$. Now, if $\omega$ is a balanced metric on $X$, then 
$$
0<\int_X i\,\varphi_1\wedge\overline{\varphi}_1\wedge\omega^2 = \int_X  d(i\,\varphi_2\wedge\omega^2) =0,
$$
which is a contradiction.

Thus, the existence of a holomorphic contact structure on a compact complex manifold  does not imply the page-$1$-$\partial\bar\partial$-property or the existence of balanced metrics. 
\end{Obs}

The examples (1) and (4) are nilmanifolds endowed with a left invariant complex structure. In the following result we prove that they are the unique complex $3$-dimensional nilmanifolds carrying a holomorphic contact structure.

\begin{Prop}\label{Prop:uniqueness-Lie-examples-dim3}
Let $X=G/\Lambda$ be the quotient by a lattice of a nilpotent Lie group of real dimension $6$ endowed with a left invariant complex structure. Suppose that $X$ has a holomorphic contact structure. Then, $X$ is the Iwasawa manifold $(1)$ or the non complex parallelisable nilmanifold $(4)$.
\end{Prop}

\noindent {\it Proof.} 
Let $\mathfrak{g}$ be the Lie algebra of the Lie group $G$, and denote by $\bigwedge^{p,\,q} \mathfrak{g}^*$ the space of left invariant forms $\widetilde\alpha$ on $G$ of bidegree $(p,\,q)$ with respect to the left invariant complex structure. There is a natural map 
$\bigwedge^{p,\,q} \mathfrak{g}^* \longrightarrow C^\infty_{p,\,q}(X,\,\C)$, sending $\widetilde\alpha$ to the corresponding induced $(p,q)$-form $\alpha$ on the quotient $X$.

In the complex dimension $3$, by the results in [Rol09] and [FRR19], for any such complex nilmanifold $X$ one has that the previous natural map 
induces an isomorphism in every 
Dolbeault cohomology group, in particular, 
$$
H^{1,\,0}_{\bar\partial}(\mathfrak{g}^*)= \{\widetilde\alpha\in \bigwedge\!^{1,\,0} \mathfrak{g}^* \mid \delbar\widetilde\alpha=0 \} \simeq H^{1,\,0}_{\bar\partial}(X,\,\C)= 
\{\alpha\in C^\infty_{1,\,0}(X,\,\C) \mid \delbar\alpha=0 \}.
$$
This implies that any holomorphic contact form $\gamma$ on $X$ is necessarily left invariant, so the problem reduces to study the existence of such contact forms on the $6$-dimensional nilpotent Lie algebra. 

Let $J$ be a complex structure on a nilpotent Lie algebra $\frg$ of real dimension 6. According to [Uga07, Proposition 2], we have: 
\begin{enumerate}
\item[{\it (a)}] If $J$ is non-nilpotent, then there is a
basis $\{\varphi_j\}_{j=1}^3$ for $\bigwedge\!^{1,\,0} \mathfrak{g}^*$ such that
$$
\left\{
\begin{array}{lcl}
d\varphi_1 \!\! & = &\!\! 0,\\[1pt]
d\varphi_2 \!\! & = &\!\! E\, \varphi_1\wedge\varphi_3 + \varphi_1\wedge\overline{\varphi}_3 \, ,\\[2pt]
d\varphi_3 \!\! & = &\!\! A\, \varphi_1\wedge\overline{\varphi}_1 + ib\,
\varphi_1\wedge\overline{\varphi}_2 - ib\bar{E}\, \varphi_2\wedge\overline{\varphi}_1 ,
\end{array}
\right.
$$
where $A,E\in \C$ with $|E|=1$, and $b\in\R-\{0\}$; 
\item[{\it (b)}] If $J$ is nilpotent, then there is a basis
$\{\varphi^j\}_{j=1}^3$ for $\bigwedge\!^{1,\,0} \mathfrak{g}^*$ satisfying
$$
\left\{
\begin{array}{lcl}
d\varphi_1 \!\! & = &\!\! 0,\\[1pt]
d\varphi_2 \!\! & = &\!\! \epsilon\, \varphi_1\wedge\overline{\varphi}_1 \, ,\\[2pt]
d\varphi_3 \!\! & = &\!\! \rho\, \varphi_1\wedge\varphi_2 + (1-\epsilon)A\,
\varphi_1\wedge\overline{\varphi}_1 + B\, \varphi_1\wedge\overline{\varphi}_2 + C\, \varphi_2\wedge\overline{\varphi}_1
+ (1-\epsilon)D\, \varphi_2\wedge\overline{\varphi}_2,
\end{array}
\right.
$$
where $A,B,C,D\in \C$, and $\epsilon,\rho \in \{0,1\}$.
\end{enumerate}

Now it is clear that in the case {\it (a)} any $(1,0)$-form $\gamma$ satisfying $\delbar\gamma=0$ is a multiple of $\varphi_1$, so $\gamma$ is closed and it does not define any holomorphic contact structure.

In the case {\it (b)}, one gets a {\it non-closed} holomorphic $(1,0)$-form only in the following two cases: $\rho=1$ and $\epsilon=A=B=C=D=0$, or $\epsilon=\rho=1$ and $B=C=0$. 
The former is the Iwasawa manifold, whereas the latter is precisely the example (4).
\hfill$\Box$

\begin{Cor}\label{Cor:p-contact-no-open}
For compact complex manifolds, the existence of holomorphic contact structures is {\bf not} a {\bf deformation open} property. 
\end{Cor}

\noindent {\it Proof.} 
This is a consequence of Proposition \ref{Prop:uniqueness-Lie-examples-dim3} applied to the Kuranishi family $(X_t)_{t}$ of small deformations of the Iwasawa manifold. In [Nak75], Nakamura divided the small deformations into three classes, namely (i), (ii) and (iii). 
All the small deformations of the Iwasawa manifold consist of left invariant complex structures (on the real nilmanifold underlying the Iwasawa manifold) that belong to the case {\it (b)} with $\epsilon=0$ in the proof of the previous proposition. Moreover, $X_t$ is complex parallelisable if and only if it belongs to the class (i). Thus, for the deformation classes (ii) and (iii), not all the coefficients $A,B,C,D$ vanish. Therefore, there do not exist holomorphic contact structures for $X_t$ in the classes (ii) and (iii).
\hfill$\Box$

 \vspace{1ex}

The examples (2) and (5) are complex $3$-dimensional solvmanifolds endowed with a left invariant complex structure such that their canonical bundle are trivialized by a left invariant holomorphic $(3,0)$-form. In the following result we prove that they are the unique 
which carry a {\it left invariant} holomorphic contact structure.

\begin{Prop}\label{Prop:uniqueness-Lie-examples-dim3-solv}
Let $X=G/\Lambda$ be the quotient by a lattice of a solvable (non nilpotent) Lie group of real dimension $6$ endowed with a left invariant complex structure with non-zero closed $(3,0)$-form. Suppose that $X$ has a left invariant holomorphic contact structure. Then, $X$ is the Nakamura manifold $(2)$ or the non complex parallelisable solvmanifold $(5)$.
\end{Prop}

\noindent {\it Proof.} 
Let $\mathfrak{g}$ be the solvable Lie algebra of the Lie group $G$, and let $\bigwedge^{p,\,q} \mathfrak{g}^*$ be the space of left invariant forms on $G$ of bidegree $(p,\,q)$ with respect to the left invariant complex structure $J$. Since $\mathfrak{g}$ is not nilpotent, one can only ensure an injection 
from $H^{1,\,0}_{\bar\partial}(\mathfrak{g}^*)$ into $H^{1,\,0}_{\bar\partial}(X,\,\C)$, so we are restricted to study the existence of holomorphic contact structures of left invariant type.

The pairs $(\frg,J)$ are classified in [FOU15] and [Tol25], and they are described in terms of their complex equations in [FOU15] and in [OU23-25]. 
If $J$ is a complex structure on $\frg$ with (non-zero) closed $(3,0)$-form, then the pair $(\frg,J)$ has one of the following equations, where $\{\omega^1,\omega^2,\omega^3\}$ always denotes a $(1,0)$-basis:
\begin{enumerate}
\item[{\it (a)}] 
\ 
$d\omega^1=A\, \omega^{1}\wedge (\omega^{3}+\omega^{\bar{3}}),\ \ 
d\omega^2=-A\, \omega^{2}\wedge (\omega^{3}+\omega^{\bar{3}}),\ \ 
d\omega^3=0,$

where $A=\cos \theta+i\sin \theta$, $\theta\in[0,\pi)$.
(These equations correspond to complex structures on two solvable Lie algebras, named $\frg_1$ and $\frg_2^{\alpha}$, where $\alpha=|\cos \theta/\sin \theta|\geq 0$ when $\theta\not=0$;
see [FOU15, Lemma 3.2]  for details).
\item[{\it (b)}] 
\ 
$d\omega^1=0,\ \ 
d\omega^2=-\frac12 \omega^{13} -(\frac12+xi)\omega^{1\bar{3}}+xi\,\omega^{3\bar{1}},\ \ 
d\omega^3=\frac12 \omega^{12} +(\frac12-\frac{i}{4x})\omega^{1\bar{2}}
+\frac{i}{4x}\,\omega^{2\bar{1}},
$

where $x\in \mathbb{R}^+$.
(These equations correspond to complex structures on the solvable Lie algebra named $\frg_3$ in [FOU15, Prop. 3.4]).
\item[{\it (c)}] 
\ 
$d\omega^1=A\,\omega^{1}\wedge (\omega^{3}+\omega^{\bar{3}}),\ \ 
d\omega^2=-A\,\omega^{2}\wedge (\omega^{3}+\omega^{\bar{3}}),\ \ 
d\omega^3=G_{11}\,\omega^{1\bar{1}}+G_{12}\,\omega^{1\bar{2}}+
\overline{G}_{12}\,\omega^{2\bar{1}}+G_{22}\,\omega^{2\bar{2}},$

where $A,G_{12}\in \mathbb{C}$ and $G_{11},G_{22}\in \mathbb{R}$, with $|A|=1$ and $(G_{11},G_{12},G_{22})\not=(0,0,0)$. (These equations correspond to complex structures on four solvable Lie algebras, named $\frg_4$, $\frg_5$, $\frg_6$ and $\frg_7$;
see [FOU15, Lemma 3.5]  for details).
\item[{\it (d)}]  
Let $\frg_8$ be the real solvable Lie algebra underlying the Nakamura manifold given in (2). By [FOU15, Prop. 3.7], all the complex structures on $\frg_8$ with closed  $(3,0)$-form are given by the following equations:

\vspace{1ex}

$(\frg_8,J) \colon \
d\omega^1= 2i\,\omega^{13} +\omega^{3\bar{3}},\ \,
d\omega^2= -2i\,\omega^{23},\ \,
d\omega^3=0;$

\vspace{1ex}

$(\frg_8,J') \colon \
d\omega^1=2i\,\omega^{13}+\omega^{3\bar{3}},\ \,
d\omega^2= -2i\,\omega^{23}+\omega^{3\bar{3}},\ \,
d\omega^3=0;$

\vspace{1ex}

$(\frg_8,J_A) \colon \ 
d\omega^1=-(A-i)\omega^{13}-(A+i)\omega^{1\bar{3}},\ \ 
d\omega^2=(A-i)\omega^{23}+(A+i)\omega^{2\bar{3}},\ \ 
d\omega^3=0,$

\hskip1.8cm where $A\in \mathbb{C}$ with ${\mathfrak I}{\frak m}\, A\not=0$.
%
\item[{\it (e)}] 
Finally, we have the equations

\vspace{.5ex}

$d\omega^1\!=\!-\omega^{3\bar{3}},\ \ 
d\omega^2\!=\frac{i}{2} \omega^{12}+\frac12 \omega^{1\bar{3}}-\frac{i}{2}\omega^{2\bar{1}},\ \ 
d\omega^3\!=\!-\frac{i}{2} \omega^{13}+\frac{i}{2}\omega^{3\bar{1}},
$

\vspace{.5ex}

\noindent 
which correspond to the (unique up to isomorphism) complex structure on the solvable Lie algebra named $\frg_9$ in [FOU15, Prop. 3.9], 
and the complex equations

\vspace{1.0ex}

$d\omega^1\!=\! \omega^{13} - \omega^{1\bar{3}} + \omega^{3\bar{2}},\ \ 
d\omega^2\!= - \omega^{23} + \omega^{2\bar{3}},\ \ 
d\omega^3\!=\! 0,
$

\vspace{.5ex}

\noindent 
which correspond to the complex structure on the solvable Lie algebra named $\frg_{10}$ in [Tol25, Lemma 3.2 and Remark 3.3], that was missing in [FOU15]. 
It is proved in [OU23-25, Prop. 1] that the complex structure is unique up to isomorphism. 
\end{enumerate}

 \vspace{1ex}

Now, it is straightforward from the equations {\it (a)}\,--{\it (e)} above that in order to get a non-closed holomorphic 1-form the only possibilities  are $(\frg_8,J')$ and $(\frg_8,J_A)$ with $A=-i$. The complex structure $J_{-i}$ is precisely the one given in (2) and $J'$ is the one given in (5).
\hfill$\Box$

\begin{Obs}\label{Obs:non-solvmanifolds-dim3} 
Concerning the example $(3)$, i.e. the non-solvable case, it was recently proved in [OU23-25] that the Lie algebra underlying $SL(2,\,\C)$ is the only unimodular non-solvable Lie algebra of real dimension $6$ admitting a complex structure with a non-zero closed $(3,0)$-form.
\end{Obs}

\section{Lie derivatives with respect to vector-valued forms}\label{section:Lie} Let $X$ be a compact complex manifold with $\mbox{dim}_\C X = n=2p+1$. Suppose there exists a {\it $p$-contact structure} $\Gamma$ on $X$.

Since $\Gamma\wedge\partial\Gamma$ is a non-vanishing holomorphic $(n,\,0)$-form on $X$, it defines a non-vanishing global holomorphic section of the canonical line bundle $K_X$. Consequently, $K_X$ is {\it trivial}, so $X$ is a Calabi-Yau manifold. After possibly multiplying $\Gamma$ by a constant, we may assume that \begin{eqnarray*}\int\limits_X dV_\Gamma =1,\end{eqnarray*} where \begin{eqnarray*}dV_\Gamma:=(i^{p^2}\,\Gamma\wedge\overline\Gamma)\wedge(i^{(p+1)^2}\,\partial\Gamma\wedge\bar\partial\overline\Gamma) = i^{n^2}\,(\Gamma\wedge\partial\Gamma)\wedge(\overline\Gamma\wedge\bar\partial\overline\Gamma)\end{eqnarray*} Having enforced this normalisation, we will call the non-vanishing holomorphic $(n,\,0)$-form $u_\Gamma:=\Gamma\wedge\partial\Gamma$ the {\it Calabi-Yau form} induced by $\Gamma$.


\subsection{Background}\label{subsection:background_C-Y}

Since $T^{1,\,0}X$ is a holomorphic vector bundle, it has a canonical $\bar\partial$-operator defining its holomorphic structure, the natural extension of the $\bar\partial$-operator acting on the scalar-valued forms on $X$ and induced by the complex structure of $X$. Meanwhile, the contraction operation $\xi\lrcorner\cdot\,:C^\infty_{p,\,q}(X,\,\C)\longrightarrow C^\infty_{p-1,\,q}(X,\,\C)$ by any vector field $\xi\in C^\infty(X,\,T^{1,\,0}X)$ applies pointwise to $\C$-valued forms of any bidegree $(p,\,q)$. The general formula \begin{eqnarray}\label{eqn:xi_contraction-Leibniz}\bar\partial(\xi\lrcorner u) = (\bar\partial\xi)\lrcorner u - \xi\lrcorner(\bar\partial u)\end{eqnarray} holds for any $u\in C^\infty_{p,\,q}(X,\,\C)$ -- see e.g. [Pop19, Lemma 4.3].

On the other hand, for any $\theta\in C^\infty_{0,\,1}(X,\,T^{1,\,0}X)$, the pointwise operation $\theta\lrcorner\cdot\,:C^\infty_{p,\,q}(X,\,\C)\longrightarrow C^\infty_{p-1,\,q+1}(X,\,\C)$ combines the contraction of any form $u\in C^\infty_{p,\,q}(X,\,\C)$ by the $T^{1,\,0}X$-part of $\theta$ with the multiplication by the $(0,\,1)$-form part of $\theta$. The general formula \begin{eqnarray}\label{eqn:theta_contraction-Leibniz}\bar\partial(\theta\lrcorner u) = (\bar\partial\theta)\lrcorner u + \theta\lrcorner(\bar\partial u)\end{eqnarray} holds for any $u\in C^\infty_{p,\,q}(X,\,\C)$ -- see e.g. [Pop19, Lemma 4.3].

\vspace{2ex}

Recall that the Lie bracket between two elements $\varphi\in C^\infty_{0,\,p}(X,\,T^{1,\,0}X)$ and $\psi\in C^\infty_{0,\,q}(X,\,T^{1,\,0}X)$, whose expressions in local holomorphic coordinates are \begin{eqnarray*}\varphi = \sum\limits_{\lambda=1}^n\varphi^\lambda\,\frac{\partial}{\partial z_\lambda}  \hspace{3ex} \mbox{and} \hspace{3ex} \psi = \sum\limits_{\lambda=1}^n\psi^\lambda\,\frac{\partial}{\partial z_\lambda},\end{eqnarray*} where the $\varphi^\lambda$'s, resp. the $\psi^\lambda$'s, are $\C$-valued $(0,\,p)$-forms, resp. $\C$-valued $(0,\,q)$-forms, is defined by \begin{eqnarray}\label{eqn:bracket_0pq-vector_def}[\varphi,\,\psi]:=\sum\limits_{\lambda,\,\mu=1}^n\bigg(\varphi^\mu\wedge\frac{\partial\psi^\lambda}{\partial z_\mu} - (-1)^{pq}\,\psi^\mu\wedge\frac{\partial\varphi^\lambda}{\partial z_\mu}\bigg)\,\frac{\partial}{\partial z_\lambda}.\end{eqnarray} Thus, $[\varphi,\,\psi]\in C^\infty_{0,\,p+q}(X,\,T^{1,\,0}X)$. This amounts to \begin{eqnarray}\label{eqn:bracket_0pq-vector_obs}[\varphi,\,\psi] = \varphi(\psi) - (-1)^{pq}\,\psi(\varphi).\end{eqnarray}

The basic properties of this Lie bracket are summed up in the following

\begin{Lem}\label{Lem:bracket_prop} ([Kod86, $\S5$, $(5.87)--(5.89)$) (i)\, $[\varphi,\,\psi]$ is independent of the choice of local coordinates $(z_1,\dots  z_n)$.

  \vspace{1ex}

  (ii)\, For all $p,q,r\in\{0,\dots ,n\}$ and all vector-valued forms $\varphi\in C^\infty_{0,\,p}(X,\,T^{1,\,0}X)$, $\psi\in C^\infty_{0,\,q}(X,\,T^{1,\,0}X)$, $\tau\in C^\infty_{0,\,r}(X,\,T^{1,\,0}X)$, the following identities hold: \\

  \vspace{1ex}

  (a)\, $[\varphi,\,\psi] = -(-1)^{pq}\,[\psi,\,\varphi]$; \hspace{3ex} (anti-commutation) \\

  (b)\, $\bar\partial[\varphi,\,\psi] = [\bar\partial\varphi,\,\psi] + (-1)^p\,[\varphi,\,\bar\partial\psi]$; \hspace{3ex} (Leibniz rule)\\

  (c)\, $(-1)^{pr}\,[[\varphi,\,\psi],\,\tau] + (-1)^{qp}\,[[\psi,\,\tau],\,\varphi] + (-1)^{rq}\,[[\tau,\,\varphi],\,\psi] = 0$. \hspace{3ex} (Jacobi identity)

\end{Lem}

Note that the above property (b) implies that $\bar\partial[\varphi,\,\psi] = 0$ whenever $\bar\partial\varphi = 0$ and $\bar\partial\psi = 0$ and that, in this case, the $\bar\partial$-cohomology class of $[\varphi,\,\psi]$ is independent of the choices of representatives of the $\bar\partial$-cohomology classes of $\varphi$ and $\psi$. Thus, the above Lie bracket operation $[\,\cdot\,,\,\cdot\,]$ makes sense at the level of $\bar\partial$-cohomology classes as well.

\vspace{2ex}

Finally, recall the {\it Calabi-Yau isomorphism} (defined for any non-vanishing holomorphic $(n,\,0)$-form $u$, in particular for our $u_\Gamma$, and for every $q\in\{0,1,\dots ,n\}$): \begin{eqnarray}\label{eqn:C_Y-isomorphism_scalar}T_\Gamma : C^\infty_{0,\,q}(X,\,T^{1,\,0}X)\stackrel{\,\cdot\lrcorner u_\Gamma\hspace{2ex}}{\longrightarrow} C^\infty_{n-1,\,q}(X,\,\C),  \hspace{3ex} T_\Gamma(\theta):= \theta\lrcorner u_\Gamma = \theta\lrcorner(\Gamma\wedge\partial\Gamma)\end{eqnarray} and the isomorphism it induces in cohomology \begin{eqnarray}\label{eqn:C_Y-isomorphism_cohom}T_{[\Gamma]} : H^{0,\,q}_{\bar\partial}(X,\,T^{1,\,0}X)\stackrel{\,\cdot\lrcorner[u_\Gamma]\hspace{2ex}}{\longrightarrow} H^{n-1,\,q}_{\bar\partial}(X,\,\C),  \hspace{3ex} T_{[\Gamma]}([\theta]_{\bar\partial}):= [\theta\lrcorner u_\Gamma]_{\bar\partial} = [\theta\lrcorner(\Gamma\wedge\partial\Gamma)]_{\bar\partial}.\end{eqnarray} Note that our hypotheses and formulae (\ref{eqn:xi_contraction-Leibniz})--(\ref{eqn:theta_contraction-Leibniz}) imply the well-definedness of $T_{[\Gamma]}$ (i.e. the fact that all the cohomology classes featuring in (\ref{eqn:C_Y-isomorphism_cohom}) are independent of the choice of representative $\theta$ of $[\theta]_{\bar\partial}$.) 

\subsection{Lie derivatives w.r.t. $T^{1,\,0}X$-valued $(0,\,1)$-forms}\label{subsection:Lie-derivatives} We digress briefly in order to introduce and analyse a technical tool that will come in handy later on. In this $\S$\ref{subsection:Lie-derivatives}, $X$ is an arbitrary $n$-dimensional compact complex manifold with $n\geq 1$ arbitrary.

In $\S4$ of [PU23], for any smooth vector field $\xi\in C^{\infty}(X,\,T^{1,\,0}X)$ of type $(1,\,0)$, the {\bf $(1,\,0)$-Lie derivative} w.r.t. $\xi$ was defined as \begin{equation}\label{eqn:Lie-deriv_1-0_def}L^{1,\,0}_{\xi}:=[\partial,\,\xi\lrcorner\cdot\,].\end{equation} This means that the linear map $L^{1,\,0}_{\xi}:C^\infty_{p,\,q}(X,\,\C)\longrightarrow C^\infty_{p,\,q}(X,\,\C)$ is given by \begin{equation}\label{eqn:Lie-deriv_1-0_def-explicit}L^{1,\,0}_{\xi}u:= \partial(\xi\lrcorner u) + \xi\lrcorner\partial u\end{equation} for any differential form $u\in C^\infty_{p,\,q}(X,\,\C)$ and any bidegree $(p,\,q)$. In particular, $L^{1,\,0}_{\xi}$ is a differential operator of order $1$ and of bidegree $(0,\,0)$.

We now introduce the analogous definition in the case where the vector field $\xi$ is replaced with a $T^{1,\,0}X$-valued $(0,\,1)$-form $\theta$.

\begin{Def}\label{Def:Lie-derivatives-mixed_def} Let $X$ be a complex manifold. For any $\theta\in C^{\infty}_{0,\,1}(X,\,T^{1,\,0}X)$, the {\bf Lie derivative} w.r.t. $\theta$ is the differential operator $L_\theta:C^\infty_{p,\,q}(X,\,\C)\longrightarrow C^\infty_{p,\,q+1}(X,\,\C)$  defined as \begin{equation}\label{eqn:Lie-derivatives-mixed_def}L_\theta:=[\partial,\,\theta\lrcorner\cdot\,].\end{equation}
\end{Def}

 This means that for any differential form $u\in C^\infty_{p,\,q}(X,\,\C)$, we have \begin{equation}\label{eqn:Lie-derivatives-mixed_def_bis}\partial(\theta\lrcorner u) = L_\theta u + \theta\lrcorner\partial u.\end{equation}

 \vspace{1ex}

 Definition (\ref{eqn:Lie-derivatives-mixed_def}) makes sense in the more general case where $\theta$ is a $T^{1,\,0}X$-valued $(0,\,s)$-form with $s$ possibly different from $1$. Then, the only thing that may change in the translation (\ref{eqn:Lie-derivatives-mixed_def_bis}) of (\ref{eqn:Lie-derivatives-mixed_def}) is a sign depending on the parity of $s$. Specifically, we get \begin{equation}\label{eqn:Lie-derivatives-mixed_def_s_bis}\partial(\theta\lrcorner u) = L_\theta u + (-1)^{s+1}\,\theta\lrcorner\partial u\end{equation} for any non-negative integer $s$, any $\theta\in C^{\infty}_{0,\,s}(X,\,T^{1,\,0}X)$ and any $u\in C^\infty_{p,\,q}(X,\,\C)$.

 \vspace{1ex}

 We start by computing $L^{1,\,0}_{\xi}$ and $L_\theta$ in local coordinates. In particular, we will see that $L^{1,\,0}_{\frac{\partial}{\partial z_k}}$ coincides with the differential operator $\partial/\partial z_k$ defined on forms by applying its analogue for functions to the coefficients of the forms written in local coordinates, namely whenever $u=\sum\limits_{|I|=p,\,|J|=q}u_{I\bar{J}}\,dz_I\wedge d\bar{z}_J$, we define \begin{eqnarray*}\frac{\partial u}{\partial z_k} = \sum\limits_{|I|=p,\,|J|=q}\frac{\partial u_{I\bar{J}}}{\partial z_k}\,dz_I\wedge d\bar{z}_J.\end{eqnarray*}

\begin{Prop}\label{Prop:Lie-derivatives_local} Let $X$ be a complex manifold with $\mbox{dim}_\C X = n$.

\vspace{1ex}

(i)\, For any $\xi\in C^{\infty}(X,\,T^{1,\,0}X)$, if $\xi = \sum\limits_{k=1}^n\xi_k\,\frac{\partial}{\partial z_k}$ is its expression in local holomorphic coordinates $(z_1,\dots , z_n)$, the following identities hold for any form $u\in C^\infty_{p,\,q}(X,\,\C)$ and any bidegree $(p,\,q)$: \begin{eqnarray}\label{eqn:Lie-deriv_1-0_local}L^{1,\,0}_{\xi}u = \sum\limits_{k=1}^n\partial\xi_k\wedge\bigg(\frac{\partial}{\partial z_k}\lrcorner u\bigg) + \sum\limits_{k=1}^n\xi_k\,\frac{\partial u}{\partial z_k} = (\partial\xi)\lrcorner u + \xi(u),\end{eqnarray} where $\partial\xi:=\sum_{k=1}^n\partial\xi_k\wedge\frac{\partial}{\partial z_k}$ is a locally defined $(1,\,0)$-form with values in $T^{1,\,0}X$, while $\xi(\cdot):=\sum_{k=1}^n\xi_k\,\frac{\partial(\cdot)}{\partial z_k}$ is a locally defined first-order differential operator acting on differential forms.

In particular, \begin{eqnarray}\label{eqn:Lie-deriv_1-0_local_particular}L^{1,\,0}_{\frac{\partial}{\partial z_k}}u = \frac{\partial u}{\partial z_k},  \hspace{6ex} k=1,\dots , n.\end{eqnarray}

\vspace{1ex}

(ii)\, For any $\theta\in C^{\infty}_{0,\,1}(X,\,T^{1,\,0}X)$, if $\theta = \sum\limits_{k,\,r=1}^n\theta^k_{\bar{r}}\,d\bar{z}_r\wedge\frac{\partial}{\partial z_k}$ is its expression in local holomorphic coordinates $(z_1,\dots , z_n)$, for any form $u\in C^\infty_{p,\,q}(X,\,\C)$ and any bidegree $(p,\,q)$ we have: \begin{eqnarray}\label{eqn:Lie-derivatives-mixed_local}L_\theta u = \sum\limits_{k,\,r=1}^n\partial\theta^k_{\bar{r}}\wedge d\bar{z}_r\wedge\bigg(\frac{\partial}{\partial z_k}\lrcorner u\bigg) - \sum\limits_{k,\,r=1}^n\theta^k_{\bar{r}}\, d\bar{z}_r\wedge\frac{\partial u}{\partial z_k} = (\partial\theta)\lrcorner u - \theta(u),\end{eqnarray} where $\partial\theta:=\sum_{k,\,r=1}^n\partial\theta^k_{\bar{r}}\wedge d\bar{z}_r\wedge\frac{\partial}{\partial z_k}$ is a locally defined $(1,\,1)$-form with values in $T^{1,\,0}X$, while $\theta(\cdot):=\sum_{k,\,r=1}^n\theta^k_{\bar{r}}\wedge d\bar{z}_r\wedge\frac{\partial}{\partial z_k}$ is a locally defined first-order differential operator acting on differential forms.

 \end{Prop}

 \noindent {\it Proof.} Let $u\in C^\infty_{p,\,q}(X,\,\C)$ be an arbitrary form and let $u=\sum\limits_{|I|=p,\,|J|=q}u_{I\bar{J}}\,dz_I\wedge d\bar{z}_J$ be its expression in local coordinates.

 (i)\, We get: $\xi\lrcorner u = \sum\limits_{|I|=p,\,|J|=q}\sum\limits_{k=1}^n\xi_k\,u_{I\bar{J}}\,\bigg(\frac{\partial}{\partial z_k}\lrcorner dz_I\bigg)\wedge d\bar{z}_J$, hence \begin{eqnarray*}\partial(\xi\lrcorner u) & = & \sum\limits_{\substack{|I|=p\\|J|=q \\ 1\leq k,\,l\leq n}} \frac{\partial\xi_k}{\partial z_l}\,dz_l\wedge\bigg(u_{I\bar{J}}\,\bigg(\frac{\partial}{\partial z_k}\lrcorner dz_I\bigg)\wedge d\bar{z}_J\bigg) + \sum\limits_{\substack{|I|=p\\|J|=q \\ 1\leq k,\,l\leq n}}\xi_k\,\frac{\partial u_{I\bar{J}}}{\partial z_l}\,dz_l\wedge\bigg(\frac{\partial}{\partial z_k}\lrcorner dz_I\bigg)\wedge d\bar{z}_J \\
   & = & \sum\limits_{1\leq k\leq n}\partial\xi_k\wedge\bigg(\frac{\partial}{\partial z_k}\lrcorner\sum\limits_{\substack{|I|=p\\|J|=q}}u_{I\bar{J}}\,dz_I\wedge d\bar{z}_J\bigg) + \sum\limits_{\substack{|I|=p\\|J|=q \\ 1\leq k,\,l\leq n}}\xi_k\,\frac{\partial u_{I\bar{J}}}{\partial z_l}\,\bigg(-\frac{\partial}{\partial z_k}\lrcorner(dz_l\wedge dz_I) + \delta_{kl}\,dz_I\bigg)\wedge d\bar{z}_J\end{eqnarray*} Formula (\ref{eqn:Lie-deriv_1-0_local}) follows from this after we notice that the last sum equals \begin{eqnarray*} & & -\sum\limits_{1\leq k\leq n}\xi_k\,\frac{\partial}{\partial z_k}\lrcorner\bigg(\sum\limits_{\substack{|I|=p\\|J|=q}}\sum\limits_{1\leq l\leq n}\frac{\partial u_{I\bar{J}}}{\partial z_l}\,dz_l\wedge dz_I\wedge d\bar{z}_J\bigg) +  \sum\limits_{\substack{|I|=p\\|J|=q \\ 1\leq k\leq n}}\xi_k\,\frac{\partial u_{I\bar{J}}}{\partial z_k}\,dz_I\wedge d\bar{z}_J = -\xi\lrcorner\partial u +  \sum\limits_{k=1}^n\xi_k\,\frac{\partial u}{\partial z_k}.\end{eqnarray*}

   In the special case when $\xi = \frac{\partial}{\partial z_k} $, all the coefficients of $\xi$ are constant, namely $\xi_l = \delta_{kl}$ for all $l=1,\dots , n$, hence $\partial\xi_l = 0$ for all $l$ and (\ref{eqn:Lie-deriv_1-0_local}) reduces to (\ref{eqn:Lie-deriv_1-0_local_particular}).

   (ii)\, We get: $\theta\lrcorner u = \sum\limits_{|I|=p,\,|J|=q}\sum\limits_{1\leq k,\,r\leq n}\theta^k_{\bar{r}}\,u_{I\bar{J}}\,d\bar{z}_r\wedge\bigg(\frac{\partial}{\partial z_k}\lrcorner dz_I\bigg)\wedge d\bar{z}_J$, hence $\partial(\theta\lrcorner u) = T_1 + T_2$ where \begin{eqnarray*}T_1 & = & \sum\limits_{\substack{|I|=p\\|J|=q \\ 1\leq k,\,l,\,r\leq n}} \frac{\partial\theta^k_{\bar{r}}}{\partial z_l}\,dz_l\wedge d\bar{z}_r\wedge\bigg(u_{I\bar{J}}\,\bigg(\frac{\partial}{\partial z_k}\lrcorner dz_I\bigg)\wedge d\bar{z}_J\bigg) = \sum\limits_{1\leq k,\,r\leq n}\partial\theta^k_{\bar{r}}\wedge d\bar{z}_r\wedge\bigg(\frac{\partial}{\partial z_k}\lrcorner\sum\limits_{\substack{|I|=p\\|J|=q}}u_{I\bar{J}}\,dz_I\wedge d\bar{z}_J\bigg),\end{eqnarray*} which amounts to $T_1 = \partial\theta\lrcorner u$, while $T_2$ equals \begin{eqnarray*} & & \sum\limits_{\substack{|I|=p\\|J|=q \\ 1\leq k,\,l,\,r\leq n}}\theta^k_{\bar{r}}\,\frac{\partial u_{I\bar{J}}}{\partial z_l}\,dz_l\wedge d\bar{z}_r\wedge\bigg(\frac{\partial}{\partial z_k}\lrcorner dz_I\bigg)\wedge d\bar{z}_J \\
     & = & \sum\limits_{\substack{|I|=p\\|J|=q \\ 1\leq k,\,l,\,r\leq n}}\theta^k_{\bar{r}}\,\frac{\partial u_{I\bar{J}}}{\partial z_l}\,\bigg(\frac{\partial}{\partial z_k}\lrcorner(dz_l\wedge d\bar{z}_r\wedge dz_I) - \delta_{kl}\,d\bar{z}_r\wedge dz_I\bigg)\wedge d\bar{z}_J  \\
     & = & \sum\limits_{1\leq k,\,r\leq n}\theta^k_{\bar{r}}\,d\bar{z}_r\wedge\bigg(\frac{\partial}{\partial z_k}\lrcorner\sum\limits_{\substack{|I|=p\\|J|=q}}\sum\limits_{1\leq l\leq n}\frac{\partial u_{I\bar{J}}}{\partial z_l}\,dz_l\wedge dz_I\wedge d\bar{z}_J\bigg) -  \sum\limits_{\substack{|I|=p\\|J|=q \\ 1\leq k,\,r\leq n}}\theta^k_{\bar{r}}\,d\bar{z}_r\wedge\frac{\partial u_{I\bar{J}}}{\partial z_k}\,dz_I\wedge d\bar{z}_J \\
     & = & \theta\lrcorner\partial u -  \sum\limits_{k,\,r=1}^n\theta^k_{\bar{r}}\,d\bar{z}_r\wedge\frac{\partial u}{\partial z_k} = \theta\lrcorner\partial u - \theta(u).\end{eqnarray*}

   This proves formula (\ref{eqn:Lie-derivatives-mixed_local}).  \hfill $\Box$

\vspace{2ex}

We now notice a few basic properties of the Lie derivative $L_\theta$ (with respect to a $T^{1,\,0}X$-valued $(0,\,1)$-form $\theta$) introduced in Definition \ref{Def:Lie-derivatives-mixed_def}. Recall that among the properties observed in Lemma 4.2. of [PU23] for the Lie derivative $L^{1,\,0}_\xi$ of (\ref{eqn:Lie-deriv_1-0_def}) are the following identities: \begin{eqnarray}\label{eqn:L_commutations} (a)\,\, [\xi\lrcorner\cdot,\,L^{1,\,0}_{\eta}] = [L^{1,\,0}_{\xi},\,\eta\lrcorner\cdot] = [\xi,\,\eta]\lrcorner\cdot \hspace{2ex} & \mbox{and} & \hspace{2ex} (b)\,\,[L^{1,\,0}_{\xi},\,L^{1,\,0}_{\eta}] = L^{1,\,0}_{[\xi,\,\eta]} \\
 \nonumber   (c)\,\,[L^{1,\,0}_{\xi},\,\bar\partial]= [\partial,\,\bar\partial\xi\lrcorner\cdot\,] \hspace{2ex} & \mbox{and} & \hspace{2ex} (d)\,\, L^{1,\,0}_{\xi}(u\wedge v) = (L^{1,\,0}_{\xi}u)\wedge v + u\wedge L^{1,\,0}_{\xi}v  \end{eqnarray} that hold for all $(1,\,0)$-vector fields $\xi,\eta\in C^{\infty}(X,\,T^{1,\,0}X)$ and all $\C$-valued differential forms $u,v$ (of any degrees). In particular, if $\xi$ is holomorphic, (c) becomes $[L^{1,\,0}_{\xi},\,\bar\partial] = 0$.

\begin{Lem}\label{Lem:Lie-derivatives-mixed_prop} Let $\theta,\psi\in C^{\infty}_{0,\,1}(X,\,T^{1,\,0}X)$. The following identities hold: \begin{eqnarray}\label{eqn:Lie-derivatives-mixed_prop_0} & & [L_\theta,\,\bar\partial] = - [(\bar\partial\theta)\lrcorner\cdot\,,\,\partial]; \hspace{2ex} \mbox{in particular},\hspace{2ex} L_\theta\bar\partial = -\bar\partial L_\theta \hspace{2ex} \mbox{when} \hspace{2ex} \bar\partial\theta = 0; \\
    \label{eqn:Lie-derivatives-mixed_prop_0-bis} & & L_\theta\partial = -\partial L_\theta; \\
    \label{eqn:Lie-derivatives-mixed_prop_1} & & L_\theta(u\wedge v) = L_\theta(u)\wedge v + (-1)^{\deg u}\,u\wedge L_\theta(v) \hspace{5ex} \mbox{for all forms}\hspace{1ex} u,v;  \\
    \label{eqn:Lie-derivatives-mixed_prop_2} & & [\theta\lrcorner\cdot\,,\,L_\psi] = -[L_\theta,\,\psi\lrcorner\cdot\,] = [\theta,\,\psi]\lrcorner\cdot\, ; \\
   \label{eqn:Lie-derivatives-mixed_prop_3} & & [L_\theta\,,\,L_\psi] = L_{[\theta,\,\psi]}.\end{eqnarray}

\vspace{1ex}

\end{Lem}

\noindent {\it Proof.} $\bullet$ From (\ref{eqn:Lie-derivatives-mixed_def}), we get: $[L_\theta,\,\bar\partial] = [[\partial,\,\theta\lrcorner\cdot\,],\,\bar\partial]$. Meanwhile, the Jacobi identity yields: \begin{eqnarray*}-\bigg[[\partial,\,\theta\lrcorner\cdot\,],\,\bar\partial\bigg] + \bigg[[\theta\lrcorner\cdot\,,\,\bar\partial],\,\partial\bigg] + \bigg[[\bar\partial,\,\partial],\,\theta\lrcorner\cdot\,\bigg] = 0.\end{eqnarray*} Now, we see that $[\bar\partial,\,\partial] = \bar\partial\partial + \partial\bar\partial = 0$, so the last term above vanishes, while $[\theta\lrcorner\cdot\,,\,\bar\partial] = -(\bar\partial\theta)\lrcorner\cdot$. This proves (\ref{eqn:Lie-derivatives-mixed_prop_0}).

\vspace{1ex}

$\bullet$ The Jacobi identity yields: \begin{eqnarray*}-\bigg[[\partial,\,\theta\lrcorner\cdot\,],\,\partial\bigg] + \bigg[[\theta\lrcorner\cdot\,,\,\partial],\,\partial\bigg] + \bigg[[\partial,\,\partial],\,\theta\lrcorner\cdot\,\bigg] = 0.\end{eqnarray*} Now, $[\partial,\,\theta\lrcorner\cdot\,] =  L_\theta$, $[\theta\lrcorner\cdot\,,\,\partial] = -[\partial,\,\theta\lrcorner\cdot\,] = -L_\theta$, while $[\partial,\,\partial] = 2\partial^2 = 0$. Thus, the Jacobi identity translates to $-2\,[L_\theta,\,\partial] = 0$, which is equivalent to (\ref{eqn:Lie-derivatives-mixed_prop_0-bis}).

\vspace{1ex}

$\bullet$ Let $k$ be the degree of $u$. Using the definition of $L_\theta$, we get: \begin{eqnarray*}L_\theta(u\wedge v) & = & \partial\bigg(\theta\lrcorner(u\wedge v)\bigg) - \theta\lrcorner\partial(u\wedge v) = \partial\bigg((\theta\lrcorner u)\wedge v + u\wedge(\theta\lrcorner v)\bigg) - \theta\lrcorner\bigg(\partial u\wedge v + (-1)^k\,u\wedge\partial v\bigg) \\
    & = & \bigg(\partial(\theta\lrcorner u)\wedge v - (\theta\lrcorner\partial u)\wedge v\bigg) + (-1)^k\,\bigg((\theta\lrcorner u)\wedge \partial v - (\theta\lrcorner u)\wedge \partial v\bigg) \\
    & + & \bigg(\partial u\wedge(\theta\lrcorner v) - \partial u\wedge(\theta\lrcorner v)\bigg) + (-1)^k\,\bigg(u\wedge\partial(\theta\lrcorner v) - u\wedge(\theta\lrcorner\partial v)\bigg) \\
  & = & L_\theta u\wedge v + (-1)^k\,u\wedge L_\theta v.\end{eqnarray*} This proves (\ref{eqn:Lie-derivatives-mixed_prop_1}).

$\bullet$ On the other hand, the Jacobi identity yields: \begin{eqnarray*}[\theta\lrcorner\cdot\, [\partial,\,\psi\lrcorner\cdot\,]] + [\partial,\,[\psi\lrcorner\cdot\,,\,\theta\lrcorner\cdot\,]] + [\psi\lrcorner\cdot\,,[\theta\lrcorner\cdot\,,\partial]] = 0.\end{eqnarray*} Now, $[\psi\lrcorner\cdot\,,\,\theta\lrcorner\cdot\,]  = \psi\lrcorner(\theta\lrcorner\cdot\,) - \theta\lrcorner(\psi\lrcorner\cdot\,) = 0$, while $[\theta\lrcorner\cdot\,,\partial] = -[\partial,\,\theta\lrcorner\cdot\,] = -L_\theta$. Thus, the above Jacobi identity amounts to $[\theta\lrcorner\cdot\,L_\psi] - [\psi\lrcorner\cdot\,L_\theta] = 0$. This proves the first identity in (\ref{eqn:Lie-derivatives-mixed_prop_2}) since $[\psi\lrcorner\cdot\,L_\theta] = - [L_\theta,\,\psi\lrcorner\cdot\,]$.

\vspace{1ex}

$\bullet$ To prove the second identity in (\ref{eqn:Lie-derivatives-mixed_prop_2}), we will now prove the identity \begin{eqnarray}\label{eqn:Lie-derivatives-mixed_prop_2_proof}[\theta\lrcorner\cdot\,,\,L_\psi] = [\theta,\,\psi]\lrcorner\cdot\,.\end{eqnarray}

In local holomorphic coordinates $z_1,\dots , z_n$, $\theta$ and $\psi$ are of the following shapes: \begin{eqnarray*}\theta = \sum\limits_{1\leq j,\,s\leq n}\theta^j_{\bar{s}}\,d\bar{z}_s\wedge\frac{\partial}{\partial z_j}  \hspace{3ex}\mbox{and}\hspace{3ex} \psi = \sum\limits_{1\leq k,\,r\leq n}\psi^k_{\bar{r}}\,d\bar{z}_r\wedge\frac{\partial}{\partial z_k}.\end{eqnarray*} Formula (\ref{eqn:Lie-derivatives-mixed_local}) combined with (\ref{eqn:Lie-deriv_1-0_local_particular}) yields: \begin{eqnarray*}L_\psi = (\partial\psi)\lrcorner\cdot - \sum\limits_{1\leq k,\,r\leq n}\psi^k_{\bar{r}}\,d\bar{z}_r\wedge L^{1,\,0}_{\frac{\partial}{\partial z_k}}(\cdot).\end{eqnarray*}

This leads to: \begin{eqnarray}\label{eqn:Lie-derivatives-mixed_prop_2_proof_1}[\theta\lrcorner\cdot\,,\,L_\psi](u) & = & [\theta\lrcorner\cdot\,,\partial\psi\lrcorner\cdot\,](u) + \sum\limits_{1\leq k,\,r\leq n}\psi^k_{\bar{r}}\,d\bar{z}_r\wedge\bigg[L^{1,\,0}_{\frac{\partial}{\partial z_k}},\,\theta\lrcorner\cdot\,\bigg](u)\end{eqnarray}

For the first term on the r.h.s. of (\ref{eqn:Lie-derivatives-mixed_prop_2_proof_1}), we get: \begin{eqnarray*}[\theta\lrcorner\cdot\,,\partial\psi\lrcorner\cdot\,] & = & \sum\limits_{j,\,k,\,r,\,s}\bigg[\theta^j_{\bar{s}}\,d\bar{z}_s\wedge\frac{\partial}{\partial z_j}\lrcorner\cdot\,,\,\partial\psi^k_{\bar{r}}\wedge d\bar{z}_r\wedge\frac{\partial}{\partial z_k}\lrcorner\cdot\,\bigg] = \sum\limits_{j,\,k,\,r,\,s}\theta^j_{\bar{s}}\,d\bar{z}_s\wedge\bigg(\frac{\partial}{\partial z_j}\lrcorner\partial\psi^k_{\bar{r}}\bigg)\wedge d\bar{z}_r\wedge\frac{\partial}{\partial z_k}\lrcorner\cdot \\
  & + &  \sum\limits_{j,\,k,\,r,\,s}\theta^j_{\bar{s}}\,d\bar{z}_s\wedge\partial\psi^k_{\bar{r}}\wedge d\bar{z}_r\wedge\bigg(\frac{\partial}{\partial z_j}\lrcorner\frac{\partial}{\partial z_k}\lrcorner\cdot\,\bigg) + \sum\limits_{j,\,k,\,r,\,s}\theta^j_{\bar{s}}\,\partial\psi^k_{\bar{r}}\wedge d\bar{z}_r\wedge d\bar{z}_s\wedge\bigg(\frac{\partial}{\partial z_k}\lrcorner\frac{\partial}{\partial z_j}\lrcorner\cdot\,\bigg).\end{eqnarray*} We see that the last two sums cancel each other, while the last sum on the first line above equals $\theta(\psi)\lrcorner\cdot$. Thus, we have got: \begin{eqnarray}\label{eqn:Lie-derivatives-mixed_prop_2_proof_2}[\theta\lrcorner\cdot\,,\partial\psi\lrcorner\cdot\,] = \theta(\psi)\lrcorner\cdot\,.\end{eqnarray}

As for the second term on the r.h.s. of (\ref{eqn:Lie-derivatives-mixed_prop_2_proof_1}), we need to compute \begin{eqnarray*}\bigg[L^{1,\,0}_{\frac{\partial}{\partial z_k}},\,\theta\lrcorner\cdot\,\bigg](u) = L^{1,\,0}_{\frac{\partial}{\partial z_k}}(\theta\lrcorner u) - \theta\lrcorner L^{1,\,0}_{\frac{\partial}{\partial z_k}}(u).\end{eqnarray*}

We get: \begin{eqnarray*}L^{1,\,0}_{\frac{\partial}{\partial z_k}}(\theta\lrcorner u) & = & \sum\limits_{j,\,s}L^{1,\,0}_{\frac{\partial}{\partial z_k}}\bigg(\theta^j_{\bar{s}}\,\,d\bar{z}_s\wedge\bigg(\frac{\partial}{\partial z_j}\lrcorner u\bigg)\bigg) \\
  & = & \sum\limits_{j,\,s}L^{1,\,0}_{\frac{\partial}{\partial z_k}}\bigg(\theta^j_{\bar{s}}\,\,d\bar{z}_s\bigg)\wedge\bigg(\frac{\partial}{\partial z_j}\lrcorner u\bigg) + \sum\limits_{j,\,s}\bigg(\theta^j_{\bar{s}}\,\,d\bar{z}_s\bigg)\wedge L^{1,\,0}_{\frac{\partial}{\partial z_k}}\bigg(\frac{\partial}{\partial z_j}\lrcorner u\bigg),\end{eqnarray*} where we applied the Leibniz rule for $L^{1,\,0}_{\partial/\partial z_k}$ (cf. (d) of (\ref{eqn:L_commutations})) to get the last equality. Now, \begin{eqnarray*}L^{1,\,0}_{\frac{\partial}{\partial z_k}}\bigg(\theta^j_{\bar{s}}\,d\bar{z}_s\bigg) = \frac{\partial}{\partial z_k}\lrcorner\partial(\theta^j_{\bar{s}}\,d\bar{z}_s) + \partial\bigg(\frac{\partial}{\partial z_k}\lrcorner(\theta^j_{\bar{s}}\,d\bar{z}_s)\bigg) = \frac{\partial}{\partial z_k}\lrcorner(\partial\theta^j_{\bar{s}}\wedge d\bar{z}_s) = \bigg(\frac{\partial}{\partial z_k}\lrcorner\partial\theta^j_{\bar{s}}\bigg)\wedge d\bar{z}_s.\end{eqnarray*} Since the last term equals $(\partial\theta^j_{\bar{s}}/\partial z_k)\wedge d\bar{z}_s$, we conclude that \begin{eqnarray*}L^{1,\,0}_{\frac{\partial}{\partial z_k}}(\theta\lrcorner u) = \sum\limits_{j,\,s}\frac{\partial\theta^j_{\bar{s}}}{\partial z_k}\,\,d\bar{z}_s\wedge\bigg(\frac{\partial}{\partial z_j}\lrcorner u\bigg) + \sum\limits_{j,\,s}\bigg(\theta^j_{\bar{s}}\,\,d\bar{z}_s\bigg)\wedge L^{1,\,0}_{\frac{\partial}{\partial z_k}}\bigg(\frac{\partial}{\partial z_j}\lrcorner u\bigg).\end{eqnarray*}

Using these computations and (\ref{eqn:Lie-derivatives-mixed_prop_2_proof_2}), (\ref{eqn:Lie-derivatives-mixed_prop_2_proof_1}) transforms to \begin{eqnarray*}[\theta\lrcorner\cdot\,,\,L_\psi](u) & = & \theta(\psi)\lrcorner u + \sum\limits_{j,\,k,\,r,\,s}\psi^k_{\bar{r}}\,\frac{\partial\theta^j_{\bar{s}}}{\partial z_k}\,d\bar{z}_r\wedge d\bar{z}_s\wedge\bigg(\frac{\partial}{\partial z_j}\lrcorner u\bigg) \\
  & + & \sum\limits_{j,\,k,\,r,\,s}\theta^j_{\bar{s}}\,\psi^k_{\bar{r}}\,d\bar{z}_r\wedge d\bar{z}_s\wedge L^{1,\,0}_{\frac{\partial}{\partial z_k}}\bigg(\frac{\partial}{\partial z_j}\lrcorner u\bigg) - \sum\limits_{j,\,k,\,r,\,s}\theta^j_{\bar{s}}\,\psi^k_{\bar{r}}\,d\bar{z}_r\wedge d\bar{z}_s\wedge\bigg(\frac{\partial}{\partial z_j}\lrcorner L^{1,\,0}_{\frac{\partial}{\partial z_k}}(u)\bigg) \\
  & = & \theta(\psi)\lrcorner u + \psi(\theta)\lrcorner u + \sum\limits_{j,\,k,\,r,\,s}\theta^j_{\bar{s}}\,\psi^k_{\bar{r}}\,d\bar{z}_r\wedge d\bar{z}_s\wedge\bigg[L^{1,\,0}_{\frac{\partial}{\partial z_k}},\,\frac{\partial}{\partial z_j}\lrcorner\cdot\,\bigg](u) \\
  & = & [\theta,\,\psi]\lrcorner u +  \sum\limits_{j,\,k,\,r,\,s}\theta^j_{\bar{s}}\,\psi^k_{\bar{r}}\,d\bar{z}_r\wedge d\bar{z}_s\wedge\bigg[\frac{\partial}{\partial z_k},\,\frac{\partial}{\partial z_j}\bigg]\lrcorner u = [\theta,\,\psi]\lrcorner u,\end{eqnarray*} where for the last but one equality we used (\ref{eqn:bracket_0pq-vector_obs}) and (a) of (\ref{eqn:L_commutations}), while the last equality follows from $[\partial/\partial z_k,\,\partial/\partial z_j] = 0$ for all $j,\,k$. This completes the proof of (\ref{eqn:Lie-derivatives-mixed_prop_2_proof}) and implicitly the one of (\ref{eqn:Lie-derivatives-mixed_prop_2}).

\vspace{1ex}

$\bullet$ It remains to prove (\ref{eqn:Lie-derivatives-mixed_prop_3}). Using definition (\ref{eqn:Lie-derivatives-mixed_def}) of $L_\theta$, we get $[L_\theta,\,L_\psi] = [[\partial,\,\theta\lrcorner\cdot\,],\,L_\psi]$. The Jacobi identity yields: \begin{eqnarray*}-\bigg[[\partial,\,\theta\lrcorner\cdot\,],\,L_\psi\bigg] + \bigg[[\theta\lrcorner\cdot\,,\,L_\psi],\,\partial\bigg] + \bigg[[L_\psi,\,\partial],\,\theta\lrcorner\cdot\,\bigg] = 0.\end{eqnarray*} Now, $[\theta\lrcorner\cdot\,,\,L_\psi] = [\theta,\,\psi]\lrcorner\cdot$ by (\ref{eqn:Lie-derivatives-mixed_prop_2}), while $[L_\psi,\,\partial] = 0$ by (\ref{eqn:Lie-derivatives-mixed_prop_0-bis}). Thus, the above Jacobi identity translates to \begin{eqnarray}\label{eqn:Lie-derivatives-mixed_prop_3_proof}[L_\theta,\,L_\psi] =  \bigg[[\theta,\,\psi]\lrcorner\cdot,\,\partial\bigg].\end{eqnarray}

Meanwhile, when evaluated on an arbitrary $\C$-valued form $\alpha$, the last expression above reads: \begin{eqnarray*}\bigg[[\theta,\,\psi]\lrcorner\cdot,\,\partial\bigg](\alpha) = [\theta,\,\psi]\lrcorner\partial\alpha + \partial\bigg([\theta,\,\psi]\lrcorner\alpha\bigg) = [\theta,\,\psi]\lrcorner\partial\alpha + L_{[\theta,\,\psi]}(\alpha) - [\theta,\,\psi]\lrcorner\partial\alpha =  L_{[\theta,\,\psi]}(\alpha),\end{eqnarray*} where for the second equality above we used definition (\ref{eqn:Lie-derivatives-mixed_def_s_bis}) for $[\theta,\,\psi]\in C^{\infty}_{0,\,2}(X,\,T^{1,\,0}X)$ (so, $s=2$).

Together with (\ref{eqn:Lie-derivatives-mixed_prop_3_proof}), this proves (\ref{eqn:Lie-derivatives-mixed_prop_3}).     \hfill $\Box$

\subsection{Applications of the Lie derivative calculus w.r.t. $T^{1,\,0}X$-valued forms}\label{subsection:applications_Lie}

As a first application of the differential calculus developed in $\S$\ref{subsection:Lie-derivatives}, we give an analogue of the standard Cartan formula expressing $(d\Omega)(\xi_0,\dots ,\xi_p)$, where $\Omega$ is a scalar-valued form and $\xi_0,\dots ,\xi_p$ are vector fields, when these vector fields are replaced by vector-valued forms. Note that the latter equality in (\ref{eqn:Cartan-formula}) can be rewritten as: \begin{eqnarray*}\xi_p\lrcorner\dots\lrcorner\xi_0\lrcorner d\Gamma & = & \sum\limits_{j=0}^p(-1)^j\,L^{1,\,0}_{\xi_j}\bigg(\xi_p\lrcorner\dots\widehat{\xi_j\lrcorner}\dots\lrcorner\xi_0\lrcorner\Gamma\bigg) \\
  & + & \sum\limits_{0\leq j<k\leq p}(-1)^{j+k}\,\xi_p\lrcorner\dots\widehat{\xi_k\lrcorner}\dots\lrcorner\widehat{\xi_j\lrcorner}\dots\lrcorner\xi_0\lrcorner[\xi_j,\,\xi_k]\lrcorner\Gamma.\end{eqnarray*}

The analogue of this identity for  vector-valued forms is contained in

\begin{Prop}\label{Prop:Cartan-formula_forms} Let $X$ be an $n$-dimensional complex manifold, let $(p,\,q)$ be an arbitrary bidegree and let $\Omega\in C^\infty_{p,\,q}(X,\,\C)$. Then, for any $\theta_0,\theta_1\in C^\infty_{0,\,1}(X,\,T^{1,\,0}X)$, the following identity holds: \begin{eqnarray}\label{eqn:Cartan-formula_forms}\theta_1\lrcorner(\theta_0\lrcorner\partial\Omega) & = & -L_{\theta_0}(\theta_1\lrcorner\Omega) - L_{\theta_1}(\theta_0\lrcorner\Omega) - [\theta_0,\,\theta_1]\lrcorner\Omega + \partial\bigg(\theta_1\lrcorner(\theta_0\lrcorner\Omega)\bigg).\end{eqnarray}

\end{Prop}  

\noindent {\it Proof.} Definition (\ref{eqn:Lie-derivatives-mixed_def}) for $L_{\theta_1}$ yields: \begin{eqnarray}\label{eqn:Cartan-formula_forms_proof_1}\partial\bigg(\theta_1\lrcorner(\theta_0\lrcorner\Omega)\bigg) = L_{\theta_1}(\theta_0\lrcorner\Omega) + \theta_1\lrcorner\partial(\theta_0\lrcorner\Omega).\end{eqnarray}

Meanwhile, the same definition (\ref{eqn:Lie-derivatives-mixed_def}) for $L_{\theta_0}$ yields: $\partial(\theta_0\lrcorner\Omega) = L_{\theta_0}(\Omega) + \theta_0\lrcorner\partial\Omega$. Taking the contraction by $\theta_1$, we get \begin{eqnarray}\label{eqn:Cartan-formula_forms_proof_2}\theta_1\lrcorner\partial(\theta_0\lrcorner\Omega) = \theta_1\lrcorner(\theta_0\lrcorner\partial\Omega) + L_{\theta_0}(\theta_1\lrcorner\Omega) + [\theta_0,\,\theta_1]\lrcorner\Omega,\end{eqnarray} where we have used property (\ref{eqn:Lie-derivatives-mixed_prop_2}) to deduce that $\theta_1\lrcorner L_{\theta_0}(\Omega) = L_{\theta_0}(\theta_1\lrcorner\Omega) + [\theta_0,\,\theta_1]\lrcorner\Omega$.

Putting together (\ref{eqn:Cartan-formula_forms_proof_1}) and (\ref{eqn:Cartan-formula_forms_proof_2}), we get (\ref{eqn:Cartan-formula_forms}).  \hfill $\Box$

\begin{Cor}\label{Cor:bracket-two-forms} Let $X$ be a complex manifold with $\mbox{dim}_\C X = n =2p+1$. Suppose that $X$ carries a holomorphic $p$-contact structure $\Gamma\in C^\infty_{p,\,0}(X,\,\C)$.

\vspace{1ex}

(i)\, If $\theta_0,\theta_1\in C^\infty_{0,\,1}(X,\,T^{1,\,0}X)$ are such that $\theta_0\lrcorner\partial\Gamma = \theta_1\lrcorner\partial\Gamma = 0$, then $[\theta_0,\,\theta_1]\lrcorner\partial\Gamma = 0$.

\vspace{1ex}

(ii)\, Let $\theta_0,\theta_1\in C^\infty_{0,\,1}(X,\,T^{1,\,0}X)$ be such that $\theta_0\lrcorner\Gamma = 0$ and $\theta_1\lrcorner\partial\Gamma = 0$. Then $[\theta_0,\,\theta_1]\lrcorner\partial\Gamma = L_{\theta_1}L_{\theta_0}\Gamma$. If, moreover, $\partial(\theta_0\lrcorner\partial\Gamma) = 0$, then $[\theta_0,\,\theta_1]\lrcorner\partial\Gamma = 0$.

\vspace{1ex}

(iii)\, Let $\theta_0,\theta_1\in C^\infty_{0,\,1}(X,\,T^{1,\,0}X)$ be such that $\theta_0\lrcorner\Gamma = 0$ and $\theta_1\lrcorner\partial\Gamma = 0$. Then $[\theta_0,\,\theta_1]\lrcorner\Gamma = -L_{\theta_0}(\theta_1\lrcorner\Gamma)$. If, moreover, $\partial(\theta_1\lrcorner\Gamma) = 0$, then $[\theta_0,\,\theta_1]\lrcorner\Gamma = 0$.

\vspace{1ex}

(iv)\, Let $\theta_0,\theta_1\in C^\infty_{0,\,1}(X,\,T^{1,\,0}X)$ be such that $\theta_0\lrcorner\Gamma = \theta_1\lrcorner\Gamma = 0$. Then, $[\theta_0,\,\theta_1]\lrcorner\Gamma = -\theta_1\lrcorner(\theta_0\lrcorner\partial\Gamma)$.

If, moreover, $\partial(\theta_0\lrcorner\partial\Gamma) = \partial(\theta_1\lrcorner\partial\Gamma) = 0$ and $\partial(\theta_1\lrcorner(\theta_0\lrcorner\partial\Gamma)) = 0$, then $[\theta_0,\,\theta_1]\lrcorner\partial\Gamma = 0$.

\end{Cor}  
  
\noindent {\it Proof.} (i)\, If we take $\Omega = \partial\Gamma$, then $\partial\Omega = 0$, $\theta_0\lrcorner\Omega = \theta_1\lrcorner\Omega = 0$ and $\theta_1\lrcorner(\theta_0\lrcorner\Omega) = 0$, so (\ref{eqn:Cartan-formula_forms}) reduces to $[\theta_0,\,\theta_1]\lrcorner\partial\Gamma = 0$, as claimed.

\vspace{1ex}

(ii)\, Taking again $\Omega = \partial\Gamma$, (\ref{eqn:Cartan-formula_forms}) reduces to $[\theta_0,\,\theta_1]\lrcorner\partial\Gamma = - L_{\theta_1}(\theta_0\lrcorner\partial\Gamma)$ under our assumptions.

Meanwhile, we have: $0 = \partial(\theta_0\lrcorner\Gamma) = L_{\theta_0}\Gamma + \theta_0\lrcorner\partial\Gamma$. In other words, we have the general fact: \begin{eqnarray}\label{eqn:L_theta_Gamma_G}\theta_0\lrcorner\partial\Gamma = - L_{\theta_0}\Gamma \hspace{5ex}\mbox{whenever}\hspace{2ex} \theta_0\lrcorner\Gamma = 0.\end{eqnarray}

Putting these two pieces of information together, we get the former contention.

On the other hand, we get: \begin{eqnarray*}0 = \partial\bigg(\theta_1\lrcorner(\theta_0\lrcorner\partial\Gamma)\bigg) = L_{\theta_1}(\theta_0\lrcorner\partial\Gamma) + \theta_1\lrcorner\partial(\theta_0\lrcorner\partial\Gamma) = L_{\theta_1}(\theta_0\lrcorner\partial\Gamma) = -L_{\theta_1}L_{\theta_0}\Gamma,\end{eqnarray*} where the first equality follows from $\theta_1\lrcorner(\theta_0\lrcorner\partial\Gamma) = \theta_0\lrcorner(\theta_1\lrcorner\partial\Gamma) = 0$, the last but one equality follows from our extra assumption $\partial(\theta_0\lrcorner\partial\Gamma) = 0$ and the last equality follows from (\ref{eqn:L_theta_Gamma_G}). Thanks to this, the latter contention now follows from the former.

\vspace{1ex}

(iii)\, Taking $\Omega = \Gamma$, (\ref{eqn:Cartan-formula_forms}) reduces to $[\theta_0,\,\theta_1]\lrcorner\Gamma = - L_{\theta_0}(\theta_1\lrcorner\Gamma)$ under our assumptions. This proves the former contention.

Meanwhile, we have: $0 = \partial\bigg(\theta_0\lrcorner(\theta_1\lrcorner\Gamma)\bigg) = L_{\theta_0}(\theta_1\lrcorner\Gamma) + \theta_0\lrcorner\partial(\theta_1\lrcorner\Gamma) = L_{\theta_0}(\theta_1\lrcorner\Gamma)$, where the first equality follows from $\theta_0\lrcorner(\theta_1\lrcorner\Gamma) = \theta_1\lrcorner(\theta_0\lrcorner\Gamma) = 0$ (having used the hypothesis $\theta_0\lrcorner\Gamma = 0$) and the last equality follows from our extra assumption $\partial(\theta_1\lrcorner\Gamma) = 0$. Hence, $L_{\theta_0}(\theta_1\lrcorner\Gamma) = 0$, so the latter contention now follows from the former.

\vspace{1ex}

(iv)\, Taking $\Omega = \Gamma$, (\ref{eqn:Cartan-formula_forms}) reduces to $[\theta_0,\,\theta_1]\lrcorner\Gamma = -\theta_1\lrcorner(\theta_0\lrcorner\partial\Gamma)$ under our assumptions. This proves the former contention.

Taking $\partial$ in the above equality, we get the second equality below: \begin{eqnarray*}L_{[\theta_0,\,\theta_1]}\Gamma - [\theta_0,\,\theta_1]\lrcorner\partial\Gamma =  \partial\bigg([\theta_0,\,\theta_1]\lrcorner\Gamma\bigg) = -\partial\bigg(\theta_1\lrcorner(\theta_0\lrcorner\partial\Gamma)\bigg) = 0,\end{eqnarray*} the last equality being a consequence of our third extra assumption. From this we get the first equality below, while from (\ref{eqn:Lie-derivatives-mixed_prop_3}) we get the second: \begin{eqnarray}\label{Cor:bracket-two-forms_proof_1}[\theta_0,\,\theta_1]\lrcorner\partial\Gamma = L_{[\theta_0,\,\theta_1]}\,\Gamma = L_{\theta_0}(L_{\theta_1}\Gamma) + L_{\theta_1}(L_{\theta_0}\Gamma).\end{eqnarray}

Now, taking $\partial$ in $\theta_0\lrcorner\Gamma = 0$, we get the first equality below: \begin{eqnarray*}0 = \partial(\theta_0\lrcorner\Gamma) = L_{\theta_0}\Gamma + \theta_0\lrcorner\partial\Gamma.\end{eqnarray*} Thus, $L_{\theta_0}\Gamma = -\theta_0\lrcorner\partial\Gamma$. Similarly, we get $L_{\theta_1}\Gamma = -\theta_1\lrcorner\partial\Gamma$. Therefore, (\ref{Cor:bracket-two-forms_proof_1}) translates to: \begin{eqnarray}\label{Cor:bracket-two-forms_proof_2}\nonumber[\theta_0,\,\theta_1]\lrcorner\partial\Gamma & = & -L_{\theta_0}(\theta_1\lrcorner\partial\Gamma) - L_{\theta_1}(\theta_0\lrcorner\partial\Gamma)\\
  & = & [\theta_0,\,\theta_1]\lrcorner\partial\Gamma - \theta_1\lrcorner L_{\theta_0}(\partial\Gamma) + [\theta_1,\,\theta_0]\lrcorner\partial\Gamma - \theta_0\lrcorner L_{\theta_1}(\partial\Gamma),\end{eqnarray} where for the latter equality we used (\ref{eqn:Lie-derivatives-mixed_prop_2}).

Now, on the one hand, $[\theta_1,\,\theta_0] = [\theta_0,\,\theta_1]$. On the other hand, our first two extra assumptions yield the first equality in each of the following two pairs of equalities: \begin{eqnarray*}0 = \partial(\theta_0\lrcorner\partial\Gamma) = L_{\theta_0}(\partial\Gamma) \hspace{3ex} \mbox{and} \hspace{3ex} 0 = \partial(\theta_1\lrcorner\partial\Gamma) = L_{\theta_1}(\partial\Gamma).\end{eqnarray*}

Thus, (\ref{Cor:bracket-two-forms_proof_2}) translates to $[\theta_0,\,\theta_1]\lrcorner\partial\Gamma = 2\,[\theta_0,\,\theta_1]\lrcorner\partial\Gamma$, yielding the stated conclusion $[\theta_0,\,\theta_1]\lrcorner\partial\Gamma = 0$. \hfill $\Box$

\vspace{2ex}

Corollary \ref{Cor:bracket-two-forms} suggests introducing the following terminology that will come in handy later on.

\begin{Prop-Def}\label{Def:constantly_horizontal-vertical} Let $X$ be a complex manifold with $\mbox{dim}_\C X = n = 2p+1$. Suppose that $X$ carries a holomorphic $p$-contact structure $\Gamma$. Let $U\subset X$ be an open subset and let $q\in\{0,\dots , n\}$.

\vspace{1ex}

$(1)$\, A form $\theta\in C_{0,\,q}^\infty(U,\,T^{1,\,0}X)$ is said to be {\bf constantly horizontal} if \begin{eqnarray*}(a)\,\theta\lrcorner\Gamma = 0 \hspace{3ex}\mbox{and}\hspace{3ex} (b)\, L_\theta(\partial\Gamma) = 0.\end{eqnarray*} 

These conditions are equivalent to $\theta\lrcorner\Gamma = 0$ and $\partial(\theta\lrcorner\partial\Gamma) = 0$.

\vspace{1ex}

$(2)$\, A form $\theta\in C_{0,\,q}^\infty(U,\,T^{1,\,0}X)$, is said to be {\bf constantly vertical} if \begin{eqnarray*}(c)\,\theta\lrcorner\partial\Gamma = 0 \hspace{3ex}\mbox{and}\hspace{3ex} (d)\, L_\theta(\Gamma) = 0.\end{eqnarray*}

These conditions are equivalent to $\theta\lrcorner\partial\Gamma = 0$ and $\partial(\theta\lrcorner\Gamma) = 0$.

\end{Prop-Def}  

\noindent {\it Proof.} The equivalence stated in $(1)$ follows from the identity \begin{eqnarray*}\partial(\theta\lrcorner\partial\Gamma) = L_\theta(\partial\Gamma) + (-1)^{q+1}\,\theta\lrcorner\partial^2\Gamma = L_\theta(\partial\Gamma).\end{eqnarray*}

Similarly, the equivalence stated in $(2)$ follows from the identity \begin{eqnarray*}\partial(\theta\lrcorner\Gamma) = L_\theta(\Gamma) + (-1)^{q+1}\,\theta\lrcorner\partial\Gamma.\end{eqnarray*} \hfill $\Box$

\vspace{2ex}

By way of a justification of the terminology introduced in the above definition, we note that condition (a) in $(1)$ requires $\theta$ to be {\it horizontal} (i.e. $\theta$ to be an ${\cal F}_\Gamma$-valued $(0,\,q)$-form), while condition (b) requires $\partial\Gamma$ to be {\it constant} in the $\theta$-direction. The roles of $\Gamma$ and $\partial\Gamma$ get permuted and the term {\it vertical} is substituted for the term {\it horizontal} (or, equivalently, the sheaf ${\cal G}_\Gamma$ is substituted for the sheaf ${\cal F}_\Gamma$) in passing from $(1)$ to $(2)$.

\vspace{2ex}

With this terminology in place, Corollary \ref{Cor:bracket-two-forms} yields the following

\begin{Cor}\label{Cor:brackets_horizontal-vertical_forms} Let $X$ be a complex manifold with $\mbox{dim}_\C X = n =2p+1$. Suppose that $X$ carries a holomorphic $p$-contact structure $\Gamma\in C^\infty_{p,\,0}(X,\,\C)$. Let $\theta_0,\theta_1\in C^\infty_{0,\,1}(X,\,T^{1,\,0}X)$.

\vspace{1ex}

(i)\, If one of the forms $\theta_0$ and $\theta_1$ is {\bf constantly horizontal} and the other {\bf vertical}, then $[\theta_0,\,\theta_1]$ is {\bf vertical}.

Similarly, if one of the forms $\theta_0$ and $\theta_1$ is {\bf constantly vertical} and the other {\bf horizontal}, then $[\theta_0,\,\theta_1]$ is {\bf horizontal}.

Finally, if one of the forms $\theta_0$ and $\theta_1$ is {\bf constantly horizontal} and the other {\bf constantly vertical}, then $[\theta_0,\,\theta_1] = 0$.
  
\vspace{1ex}

(ii)\, If $\theta_0$ and $\theta_1$ are both {\bf constantly horizontal} and if either $L_{\theta_1}(\theta_0\lrcorner\partial\Gamma) = 0$ or $L_{\theta_0}(\theta_1\lrcorner\partial\Gamma) = 0$, then $[\theta_0,\,\theta_1]$ is {\bf vertical}.

\end{Cor}

\noindent {\it Proof.} The first two statements under (i) are respective rewordings of (ii) and (iii) of Corollary \ref{Cor:bracket-two-forms}. The third statement under (i) is an immediate consequence of the first two: these imply, under the hypotheses of the third statement, that $[\theta_0,\,\theta_1]\lrcorner\partial\Gamma = 0$ and $[\theta_0,\,\theta_1]\lrcorner\Gamma = 0$, hence also that  $[\theta_0,\,\theta_1]\lrcorner u_\Gamma = 0$. This last equality implies $[\theta_0,\,\theta_1] = 0$ thanks to the Calabi-Yau isomorphism (\ref{eqn:C_Y-isomorphism_scalar}).

Statement (ii) is a rewording of (iv) of Corollary \ref{Cor:bracket-two-forms}. (Note that $\theta_1\lrcorner(\theta_0\lrcorner\partial\Gamma) = \theta_0\lrcorner(\theta_1\lrcorner\partial\Gamma)$.)  \hfill $\Box$

\vspace{2ex}

Another application of the differential calculus for the Lie derivative $L_\theta$ w.r.t. a vector-valued form $\theta\in C^\infty(X,\,T^{1,\,0}X)$ is its extendability from forms to cohomology thanks to

\begin{Prop}\label{Prop:L_theta_well-defined_cohomology} Let $X$ be a complex manifold. For any $\theta\in C^\infty_{0,\,1}(X,\,T^{1,\,0}X)$ such that $\bar\partial\theta=0$, any bidegree $(p,\,q)$ and any $u\in C^\infty_{p,\,q}(X,\,\C)$ such that $\bar\partial u=0$, the following equality holds: \begin{eqnarray}\label{eqn:L_theta_well-defined_cohomology}\bigg[L_{\theta + \bar\partial\xi}(u+\bar\partial v)\bigg]_{\bar\partial} = [L_\theta u]_{\bar\partial},\end{eqnarray} for every $\xi\in C^\infty(X,\,T^{1,\,0}X)$ and every $v\in C^\infty_{p,\,q-1}(X,\,\C)$.

\end{Prop}

\noindent {\it Proof.} Definition (\ref{eqn:Lie-derivatives-mixed_def_bis}) yields the first equality below: \begin{eqnarray*}L_{\theta + \bar\partial\xi}(u+\bar\partial v) & = & \partial\bigg((\theta + \bar\partial\xi)\lrcorner(u+\bar\partial v)\bigg) - (\theta + \bar\partial\xi)\lrcorner\partial(u+\bar\partial v) \\
  & = & \bigg(\partial(\theta\lrcorner u) - \theta\lrcorner\partial u\bigg) + \bigg(\partial(\theta\lrcorner\bar\partial v) - \theta\lrcorner\partial\bar\partial v\bigg) + \partial\bar\partial\bigg(\xi\lrcorner(u + \bar\partial v)\bigg) - \bar\partial\bigg(\xi\lrcorner\partial(u+\bar\partial v)\bigg) \\
  & = & L_\theta u + L_\theta(\bar\partial v) - \bar\partial\bigg(\partial\bigg(\xi\lrcorner(u + \bar\partial v)\bigg) + \xi\lrcorner\partial(u+\bar\partial v)\bigg) \\
  & = & L_\theta u - \bar\partial(L_\theta v) - \bar\partial\bigg(L^{1,\,0}_\xi(u+\bar\partial v))\bigg),\end{eqnarray*} where for the last equality we used (\ref{eqn:Lie-derivatives-mixed_prop_0}) in the case $\bar\partial\theta = 0$, as well definition (\ref{eqn:Lie-deriv_1-0_def-explicit}) of $L^{1,\,0}_\xi$.

This proves the general formula: \begin{eqnarray*}L_{\theta + \bar\partial\xi}(u+\bar\partial v) =  L_\theta u - \bar\partial\bigg(L_\theta v + L^{1,\,0}_\xi(u+\bar\partial v))\bigg),\end{eqnarray*} which, in turn, proves (\ref{eqn:L_theta_well-defined_cohomology}).  \hfill $\Box$

\subsection{Generalised Tian-Todorov formulae}\label{subsection:generalised_Tian-Todorov}

The following classical result (cf. Lemma 3.1. in [Tia87], Lemma 1.2.4. in [Tod89]) makes implicit use of the back-and-forth between (cohomology classes of) $T^{1,\,0}X$-valued forms and scalar-valued forms made possible by the Calabi-Yau isomorphisms $T_\Gamma$ and $T_{[\Gamma]}$.

\begin{Lem}(Tian-Todorov Lemma)\label{Lem:Tian-Todorov} Let $X$ be a compact complex manifold ($n=\mbox{dim}_{\C}X$) such that $K_X$ is {\bf trivial}. Fix a non-vanishing holomorphic $(n,\,0)$-form $u$ on $X$.

  Then, for any forms $\theta_1, \theta_2\in C^{\infty}_{0,\, 1}(X,\, T^{1,\, 0}X)$ such that $\partial(\theta_1\lrcorner u)=\partial(\theta_2\lrcorner u)=0$, we have $$[\theta_1,\, \theta_2]\lrcorner u\in\mbox{Im}\,\partial.$$ More precisely, the identity \begin{equation}\label{eqn:basic-trick}[\theta_1,\, \theta_2]\lrcorner u = -\partial\bigg(\theta_1\lrcorner(\theta_2\lrcorner u)\bigg)\end{equation} holds for $\theta_1, \theta_2\in C^{\infty}_{0,\, 1}(X,\, T^{1,\, 0}X)$ whenever $\partial(\theta_1\lrcorner u) = \partial(\theta_2\lrcorner u)=0$.

\end{Lem}

\vspace{2ex}

We shall now observe that the properties of the Lie derivative w.r.t. a $T^{1,\,0}X$-valued $(0,\,1)$-form introduced and studied in $\S$\ref{subsection:Lie-derivatives} yield at once a more general form of the above classical result in which the scalar-valued forms $\theta_1\lrcorner u$ and $\theta_2\lrcorner u$ are not assumed $\partial$-closed. We need not even suppose $K_X$ to be trivial, as part (ii) (generalising the Calabi-Yau case part (i)) of the following result shows.

\begin{Lem}(generalised Tian-Todorov)\label{Lem:generalised_Tian-Todorov} Let $X$ be a compact complex manifold with $n=\mbox{dim}_{\C}X$.

\vspace{1ex}

(i)\, Suppose that $K_X$ is {\bf trivial} and fix a non-vanishing holomorphic $(n,\,0)$-form $u$ on $X$.

Then, for any forms $\theta_1, \theta_2\in C^{\infty}_{0,\, 1}(X,\, T^{1,\, 0}X)$, the following identity holds: \begin{equation}\label{eqn:generalised_basic-trick}[\theta_1,\, \theta_2]\lrcorner u = -\partial\bigg(\theta_1\lrcorner(\theta_2\lrcorner u)\bigg) + \theta_1\lrcorner\partial(\theta_2\lrcorner u) + \theta_2\lrcorner\partial(\theta_1\lrcorner u).\end{equation}

\vspace{1ex}

(ii)\, For any bidegree $(p,\,q)$, any form $\alpha\in C^\infty_{p,\,q}(X,\,\C)$ and any forms $\theta_1, \theta_2\in C^{\infty}_{0,\, 1}(X,\, T^{1,\, 0}X)$, the following identity holds: \begin{equation}\label{eqn:generalised_basic-trick_non-C-Y}[\theta_1,\, \theta_2]\lrcorner\alpha = -\partial\bigg(\theta_1\lrcorner(\theta_2\lrcorner\alpha)\bigg) + \theta_1\lrcorner L_{\theta_2}(\alpha) + \theta_2\lrcorner L_{\theta_1}(\alpha) + \theta_1\lrcorner(\theta_2\lrcorner\partial\alpha).\end{equation}

In particular, if $\partial\alpha = 0$ and $\partial(\theta_1\lrcorner\alpha) = \partial(\theta_2\lrcorner\alpha) = 0$, then \begin{equation}\label{eqn:generalised_basic-trick_non-C-Y_re}[\theta_1,\, \theta_2]\lrcorner\alpha = -\partial\bigg(\theta_1\lrcorner(\theta_2\lrcorner\alpha)\bigg).\end{equation}

\end{Lem}

\noindent {\it Proof.} Note that $L_{\theta_2}(\alpha) = \partial(\theta_2\lrcorner\alpha) - \theta_2\lrcorner\partial\alpha$ and the analogous equality with $\theta_1$ in place of $\theta_2$ hold. Thus, in the special case where $\partial\alpha = 0$ (which occurs when $K_X$ is trivial and $\alpha = u\in C^\infty_{n,\,0}(X,\,\C)$ is a Calabi-Yau form on $X$), identity (\ref{eqn:generalised_basic-trick_non-C-Y}) reduces to identity (\ref{eqn:generalised_basic-trick}). Therefore, it suffices to prove (ii).

Definition \ref{Def:Lie-derivatives-mixed_def} yields the first equality below: \begin{eqnarray*}\partial\bigg(\theta_1\lrcorner(\theta_2\lrcorner\alpha)\bigg) & = & L_{\theta_1}(\theta_2\lrcorner\alpha) + \theta_1\lrcorner\partial(\theta_2\lrcorner\alpha) = [L_{\theta_1},\,\theta_2\lrcorner\cdot\,](\alpha) + \theta_2\lrcorner L_{\theta_1}(\alpha) + \theta_1\lrcorner L_{\theta_2}(\alpha) + \theta_1\lrcorner(\theta_2\lrcorner\partial\alpha) \\
  & = & -[\theta_1,\,\theta_2]\lrcorner\alpha + \theta_2\lrcorner L_{\theta_1}(\alpha) + \theta_1\lrcorner L_{\theta_2}(\alpha) + \theta_1\lrcorner(\theta_2\lrcorner\partial\alpha),\end{eqnarray*} where the last equality follows from (\ref{eqn:Lie-derivatives-mixed_prop_2}).

This proves (\ref{eqn:generalised_basic-trick_non-C-Y}).   \hfill $\Box$

\vspace{2ex}

A consequence of this is the following partial strengthening of (iv) of Corollary \ref{Cor:bracket-two-forms}.

\begin{Cor}\label{Cor:generalised_Tian-Todorov_del-Gamma} Let $X$ be a complex manifold with $\mbox{dim}_\C X = n =2p+1$. Suppose that $X$ carries a holomorphic $p$-contact structure $\Gamma\in C^\infty_{p,\,0}(X,\,\C)$.

  For any forms $\theta_1, \theta_2\in C^{\infty}_{0,\, 1}(X,\, T^{1,\, 0}X)$ such that $\partial(\theta_1\lrcorner\partial\Gamma) = \partial(\theta_2\lrcorner\partial\Gamma) = 0$ and $\partial(\theta_1\lrcorner(\theta_2\lrcorner\partial\Gamma)) = 0$, we have \begin{equation}\label{eqn:generalised_Tian-Todorov_del-Gamma}[\theta_1,\, \theta_2]\lrcorner\partial\Gamma = 0.\end{equation} 

\end{Cor}

\noindent {\it Proof.} This follows at once from (\ref{eqn:generalised_basic-trick_non-C-Y_re}) by taking $\alpha = \partial\Gamma$ in (ii) of Lemma \ref{Lem:generalised_Tian-Todorov}.  \hfill $\Box$

\section{The sheaves ${\cal F}_\Gamma$ and ${\cal G}_\Gamma$: examples and computations}\label{section:sheaves_F-G_examples-computations}

Two subsheaves of the holomorphic tangent sheaf ${\cal O}(T^{1,\,0}X)$ that are naturally associated with the forms $\Gamma$ and $\partial\Gamma$ of a holomorphic $p$-contact structure were introduced in [KPU25] as follows.

\begin{Def}\label{Def:sheaves_F-G} Suppose $\Gamma\in C^\infty_{p,\,0}(X,\,\C)$ is a holomorphic $p$-contact structure on a compact complex manifold $X$ with $\mbox{dim}_\C X = n = 2p+1$.

\vspace{1ex}

(i)\, We let ${\cal F}_\Gamma$ be the sheaf of germs of holomorphic $(1,\,0)$-vector fields $\xi$ such that $\xi\lrcorner\Gamma=0$.

\vspace{1ex}

(ii)\, We let ${\cal G}_\Gamma$ be the sheaf of germs of holomorphic $(1,\,0)$-vector fields $\xi$ such that $\xi\lrcorner\partial\Gamma=0$.

\end{Def}

\vspace{2ex}

These are torsion-free, coherent, but not necessarily locally free, sheaves of ${\cal O}_X$-modules whose ranks are unpredictable.

To illustrate the sheaves ${\cal F}_\Gamma$ and ${\cal G}_\Gamma$, we now compute them in two classes of explicit examples of holomorphic $p$-contact manifolds $(X,\,\Gamma)$.

\vspace{1ex}

The first such class consists of two $7$-dimensional generalisations of the classical $3$-dimensional Iwasawa manifold $I^{(3)}$. They are slightly different from the generalisations given in [KPU25, Prop.~3.2].

\begin{Ex}\label{Ex:I_3_7-dim-analogues} For $\varepsilon\in\{0,1\}$, let $G_\varepsilon$ be the nilpotent $7$-dimensional complex Lie group whose complex structure is defined by the following structure equations involving a basis of holomorphic $(1,\,0)$-forms $\varphi_1,\dots , \varphi_7$:
$$
d\varphi_1 = d\varphi_2 = 0,  \ 
   d\varphi_3 = \varphi_1\wedge\varphi_2,  \ 
   d\varphi_4 = \varphi_1\wedge\varphi_3,  \ 
   d\varphi_5 = \varphi_2\wedge\varphi_3,  \ 
   d\varphi_6 = \varepsilon\, \varphi_2\wedge\varphi_5,  \ 
   d\varphi_7 = \varphi_2\wedge\varphi_6.
$$
Let $\Lambda_\varepsilon\subset G_\varepsilon$ be a co-compact lattice and consider the compact nilmanifold $X_\varepsilon=G_\varepsilon/\Lambda_\varepsilon$. On $X_\varepsilon$ we define the 
$(3,\,0)$-form 
$$
\Gamma_\varepsilon = \varphi_3\wedge(\varphi_4\wedge\varphi_5+ \varphi_5\wedge\varphi_6 + \varphi_6\wedge\varphi_7).
$$

Since all the forms $\varphi_j$ are holomorphic, $\delbar\Gamma_\varepsilon=0$. Moreover, a straightforward calculation gives 
$$
\partial\Gamma_\varepsilon= \varphi_1\wedge\varphi_2\wedge(\varphi_4\wedge\varphi_5+ \varphi_5\wedge\varphi_6 + \varphi_6\wedge\varphi_7) +\varepsilon\, \varphi_2\wedge\varphi_3\wedge\varphi_5\wedge\varphi_7,
$$
so $\Gamma_\varepsilon\wedge \partial\Gamma_\varepsilon = \varphi_1\wedge\varphi_2\wedge\varphi_3\wedge(\varphi_4\wedge\varphi_5+ \varphi_5\wedge\varphi_6 + \varphi_6\wedge\varphi_7)^2= 2\, \varphi_1\wedge\cdots\wedge\varphi_7 \not=0$. Therefore, $\Gamma_\varepsilon$ is a holomorphic 3-contact structure on the compact complex manifold $X_\varepsilon$ for every $\varepsilon\in \{0,1\}$. Moreover, one can prove that $X_\varepsilon$ does not admit any holomorphic contact structure. 

Let $\{ \xi_i \}_{i=1}^7$ be the basis of holomorphic $(1,\,0)$-vector fields on $X_\varepsilon$ dual to $\{ \varphi_i \}_{i=1}^7$. As for the sheaf ${\cal F}_{\Gamma_\varepsilon}$, it is easy to check that ${\cal F}_{\Gamma_\varepsilon}= \langle \xi_1, \xi_2 \rangle$ for $\varepsilon\in \{0,1\}$. Next we show that the sheaf ${\cal G}_{\Gamma_\varepsilon}$ depends on the value of $\varepsilon$.
A direct calculation gives: 
\begin{eqnarray*}
&& \xi_1\lrcorner\partial\Gamma_\varepsilon= \varphi_2\wedge\varphi_4\wedge\varphi_5 + \varphi_2\wedge\varphi_5\wedge\varphi_6 + \varphi_2\wedge\varphi_6\wedge\varphi_7,\\
&&\xi_2\lrcorner\partial\Gamma_\varepsilon= -\varphi_1\wedge\varphi_4\wedge\varphi_5 - \varphi_1\wedge\varphi_5\wedge\varphi_6 - \varphi_1\wedge\varphi_6\wedge\varphi_7 + \varepsilon\, \varphi_3\wedge\varphi_5\wedge\varphi_7,\\
&&\xi_3\lrcorner\partial\Gamma_\varepsilon= -\varepsilon\, \varphi_2\wedge\varphi_5\wedge\varphi_7,\\
&&\xi_4\lrcorner\partial\Gamma_\varepsilon= \varphi_1\wedge\varphi_2\wedge\varphi_5,\\
&&\xi_5\lrcorner\partial\Gamma_\varepsilon= -\varphi_1\wedge\varphi_2\wedge\varphi_4 + \varphi_1\wedge\varphi_2\wedge\varphi_6 + \varepsilon\, \varphi_2\wedge\varphi_3\wedge\varphi_7,\\
&&\xi_6\lrcorner\partial\Gamma_\varepsilon= -\varphi_1\wedge\varphi_2\wedge\varphi_5 + \varphi_1\wedge\varphi_2\wedge\varphi_7,\\
&&\xi_7\lrcorner\partial\Gamma_\varepsilon= -\varphi_1\wedge\varphi_2\wedge\varphi_6 - \varepsilon\, \varphi_2\wedge\varphi_3\wedge\varphi_5.
\end{eqnarray*}
Therefore, for $\varepsilon=0$ we have ${\cal G}_{\Gamma_0}= \langle \xi_3 \rangle$, whereas for $\varepsilon=1$ the sheaf is ${\cal G}_{\Gamma_1}= \{ 0 \}$. 

\end{Ex}

\vspace{1ex}

The second class of examples of holomorphic $p$-contact manifolds $(X,\,\Gamma)$ in which we compute the sheaves ${\cal F}_{\Gamma}$ and ${\cal G}_{\Gamma}$ is the one described in [KPU25, Theorem 3.5].

\begin{Ex}\label{Ex:hol-s-symplectic_hol-p-contact_F-G} 
Let $X$ be the $(4l+3)$-dimensional compact complex nilmanifold constructed in $(2)$ of 
Theorem 3.5 of [KPU25]. 

We first recall that the structure of the manifold $X$ is given as follows. 
Let $G$ be a nilpotent complex Lie group with $\mbox{dim}_\C G = n = 2s = 4l$ and let $\Lambda$ be a co-compact lattice in $G$. Denote by $Y=G/\Lambda$ the induced quotient compact complex $n$-dimensional manifold. 
Let $\{\varphi_1,\dots, \varphi_{4l}\}$ be any basis of holomorphic $(1,0)$-forms on $Y$, which are necessarily induced by left-invariant holomorphic forms on $G$,  and  we still denote by $\{\varphi_1,\dots , \varphi_{4l}\}$ the induced $\C$-basis of $H^{1,\,0}_{\bar\partial}(Y,\,\C)$. Consider the holomorphic $s$-symplectic structure  on $Y$ defined by the $(s,\,0)$-form \begin{eqnarray*}\Omega:=\varphi_1\wedge\dots\wedge\varphi_{2l} + \varphi_{2l+1}\wedge\dots\wedge\varphi_{4l}.\end{eqnarray*} 

Set $p:=2l+1$ and consider the $(2p+1)$-dimensional compact complex nilmanifold $X$ defined by a basis $\{\pi^\star\varphi_1,\dots , \pi^\star\varphi_{4l}, \varphi_{4l+1}, \varphi_{4l+2}, \varphi_{4l+3}\}$ of holomorphic $(1,\,0)$-forms whose first $4l$ members are the pullbacks under the natural projection $\pi:X\longrightarrow Y$ of the forms considered above and the three extra members satisfy the structure equations on $X$: 
$$
\partial\varphi_{4l+1} = \partial\varphi_{4l+2} = 0 \hspace{3ex}  \mbox{and} \hspace{3ex} \partial\varphi_{4l+3} = \varphi_{4l+1}\wedge\varphi_{4l+2} + \pi^\star\sigma,
$$
where $\sigma$ is any rational $d$-closed $(2,\,0)$-form on $Y$. 
Let us consider on $X$ the holomorphic $p$-contact structure defined by the 
$(p,\,0)$-form \begin{eqnarray*}\Gamma:=\pi^\star\Omega\wedge\varphi_{4l+3}.\end{eqnarray*}


Let $\{\xi_1,\dots, \xi_{4l+3}\}$ be the basis of holomorphic $(1,\,0)$-vector fields on $X$ dual to the basis of holomorphic $(1,\,0)$-forms 
$\{\pi^\star\varphi_1,\dots , \pi^\star\varphi_{4l}, \varphi_{4l+1}, \varphi_{4l+2}, \varphi_{4l+3}\}$ on $X$. Then:

\vspace{1ex}

(i)\, the sheaf ${\cal F}_\Gamma$ of $(X,\,\Gamma)$ is the trivial rank-two holomorphic vector bundle on $X$ generated at every point by the global holomorphic $(1,\,0)$-vector fields $\xi_{4l+1}$ and $\xi_{4l+2}$;

\vspace{1ex}

(ii)\, the sheaf ${\cal G}_\Gamma$ of $(X,\,\Gamma)$ depends on the holomorphic $s$-symplectic structure $\Omega$ of $Y$, but not on $\sigma$, in the following way:

\vspace{1ex}

(a)\, if $\partial\Omega = 0$, then ${\cal G}_\Gamma$ is the trivial holomorphic line bundle on $X$ generated at every point by the global holomorphic $(1,\,0)$-vector field $\xi_{4l+3}$;

\vspace{1ex}

(b)\, if $\partial\Omega\neq 0$, then ${\cal G}_\Gamma = 0$.

\end{Ex}

\noindent {\it Proof.}   
Recall that  $\Gamma=\pi^\star\Omega\wedge\varphi_{4l+3}$, where $\Omega=\varphi_1\wedge\dots\wedge\varphi_{2l} + \varphi_{2l+1}\wedge\dots\wedge\varphi_{4l}$.

Note that the $4l$ holomorphic $(2l\!-\!1,0)$-forms: 
\begin{equation}\label{contractions-xi-s}
\begin{array}{rl}
& \xi_i\lrcorner\pi^\star\Omega = (-1)^{i+1}\pi^\star\varphi_1 \wedge\ldots\wedge\widehat\pi^\star\varphi_i\wedge\ldots\wedge \pi^\star\varphi_{2l}, \mbox{ for } 1\leq i\leq 2l,\\[6pt] 
& \xi_j\lrcorner\pi^\star\Omega = (-1)^{j+1}\pi^\star\varphi_{2l+1} \wedge\ldots\wedge\widehat\pi^\star\varphi_j\wedge\ldots\wedge \pi^\star\varphi_{4l},  \mbox{ for } 2l+1\leq j\leq 4l,
\end{array}
\end{equation}
are linearly independent at every point of $X$.  
Therefore,  the collection of $4l+1$ forms made up of $\xi_i\lrcorner\Gamma=(\xi_i\lrcorner\pi^\star\Omega) \wedge\varphi_{4l+3}$, for $1\leq i\leq 2l$, $\xi_j\lrcorner\Gamma=(\xi_j\lrcorner\pi^\star\Omega)\wedge\varphi_{4l+3}$, for $2l+1\leq j\leq 4l$, and 
$\xi_{4l+3}\lrcorner\Gamma=\pi^\star\Omega$ is linearly independent at every point of $X$.  On the other hand, by the definition of $\Gamma$, we have $\xi_{4l+1}\lrcorner\Gamma=0=\xi_{4l+2}\lrcorner\Gamma$, so the proof of (i) is complete.

Since $\partial\Gamma =   (\pi^\star\partial\Omega)\wedge\varphi_{4l+3} + \pi^\star\Omega\wedge\varphi_{4l+1}\wedge\varphi_{4l+2} + \pi^\star(\Omega\wedge\sigma)$, we get
$$
\xi_k\lrcorner\partial\Gamma =  (\xi_k\lrcorner\pi^\star\partial\Omega)\wedge\varphi_{4l+3} + (\xi_k\lrcorner\pi^\star\Omega)\wedge\varphi_{4l+1}\wedge\varphi_{4l+2} + \xi_k\lrcorner\pi^\star(\Omega\wedge\sigma), \ \mbox{ for } 1\leq k\leq 4l.
$$
Note that both terms $\pi^\star\partial\Omega$ and $\pi^\star(\Omega\wedge\sigma)$ belong to $\bigwedge\langle \pi^\star\varphi_1,\dots , \pi^\star\varphi_{4l}\rangle$, so using again 
\eqref{contractions-xi-s} we arrive at the conclusion that the following $4l+2$ forms: 
$$
\xi_k\lrcorner\partial\Gamma \ \ (1\leq k\leq 4l),\quad  \xi_{4l+1}\lrcorner\partial\Gamma= \pi^\star\Omega\wedge\varphi_{4l+2}, \quad 
\xi_{4l+2}\lrcorner\partial\Gamma=- \pi^\star\Omega\wedge\varphi_{4l+1},
$$
are linearly independent at every point of $X$. 
Finally, notice that $\xi_{4l+3}\lrcorner\partial\Gamma=\pi^\star\partial\Omega$ is again linearly independent of the previous collection if $\partial\Omega\ne0$. 
This completes the proof of (ii). 

\hfill $\Box$

\section{Contact hyperbolicity}\label{section:contact-hyperbolicity} In this section, we introduce two notions of hyperbolicity for holomorphic $p$-contact manifolds, one for each of the two aspects (metric and based on entire maps) of the theory.

\vspace{1ex}

On the metric side of the theory, we propose the following

\begin{Def}\label{Def:p-contact-metric-hyperbolic} Let $X$ be a compact complex manifold with $\mbox{dim}_\C X = n = 2p+1$. Suppose that $X$ carries a holomorphic $p$-contact structure $\Gamma\in C^\infty_{p,\,0}(X,\,\C)$.

  We say that $(X,\,\Gamma)$ is a {\bf $p$-contact metrically hyperbolic} manifold if there exist a Hermitian metric $\omega$ and a real $C^\infty$ form $\Omega$ of bidegree $(p,\,p)$ on $X$ such that:

\vspace{1ex}

(i)\, $d\Omega = 0$ and $\Omega$ is weakly strictly positive on $X$;

\vspace{1ex}

(ii)\, $\displaystyle\frac{\omega^p}{p!} = i^{p^2}\,\Gamma\wedge\overline\Gamma + \Omega$ on $X$;

\vspace{1ex}

(iii)\, $\Omega$ is $\widetilde{d}(\mbox{bounded})$ in the usual sense that, if $\pi_X:\widetilde{X}\longrightarrow X$ is the universal covering map of $X$, there exists a $C^\infty$ $(2p-1)$-form $\widetilde\beta$ on $\widetilde{X}$ that is bounded with respect to the pullback metric $\widetilde\omega:=\pi_X^\star\omega$ and has the property $\pi_X^\star\Omega = d\widetilde\beta$ on $\widetilde{X}$.

\end{Def}

\vspace{1ex}

On the entire-map-based side of the theory, we propose the following

\begin{Def}\label{Def:p-contact-hyperbolic} Let $X$ be a compact complex manifold with $\mbox{dim}_\C X = n = 2p+1$. Suppose that $X$ carries a holomorphic $p$-contact structure $\Gamma\in C^\infty_{p,\,0}(X,\,\C)$.

  We say that $(X,\,\Gamma)$ is a {\bf $p$-contact hyperbolic} manifold if there exists no holomorphic map $f:\C^p\longrightarrow X$ satisfying the following three conditions:

  \vspace{1ex}

  (i)\, $f$ is non-degenerate at some point $x_0\in\C^p$ in the usual sense that its differential map $d_{x_0}f:\C^p\longrightarrow T_{f(x_0)}^{1,\,0}X$ at $x_0$ is injective;

  \vspace{1ex}

  (ii)\, $f$ has subexponential growth (in the sense of Definition 2.2. in [KP23]);

  \vspace{1ex}

  (iii)\, $f^\star\Gamma = 0$. (We say in this case that $f$ is {\bf horizontal}.)

\end{Def}

The notion of {\it subexponential growth} was introduced in [MP22, Definition 2.3.] in the case of maps $f:\C^{n-1}\longrightarrow X$ and repeated in [KP23, Definition 2.2.] in the general case of maps $f:\C^q\longrightarrow X$ for an arbitrary $q\in\{1,\dots , n-1\}$. For the reader's convenience, we reproduce it here
together with the following notation.

\vspace{1ex}

$\bullet$ Let $B_r$, respectively $S_r$, be the open ball, respectively the sphere, of radius $r>0$ centred at the origin in $\C^q$. The $(\omega,\,f)$-volume of $B_r$, respectively the $(\omega,\,f)$-area of $S_r$, is defined by \begin{eqnarray*}\mbox{Vol}_{\omega,\,f}(B_r):=\int\limits_{B_r}f^\star\omega_q,   \hspace{3ex}\mbox{respectively}\hspace{3ex} \mbox{A}_{\omega,\,f}(S_r):=\int\limits_{S_r}d\sigma_{\omega,\,f,\,r},\end{eqnarray*} where we put $\omega_q:=\omega^q/q!$ and $d\sigma_{\omega,\,f,\,r}$ is the area measure induced by $f^\star\omega$ on the spheres of $\C^q$. It is computed using formula (3) in [KP23] and has the explicit expression (4) in [KP23].

\vspace{1ex}

$\bullet$ Let $\tau$ be the function on $\C^q$ defined by $\tau(z)=|z|^2$, the squared Euclidean norm of any $z\in\C^q$.

\begin{Def}\label{Def:subexp_re}([KP23]) Let $(X,\,\omega)$ be a compact complex Hermitian manifold with $\mbox{dim}_\C X =n\geq 2$. Let $q\in\{1,\dots , n-1\}$ and let $f:\C^q\longrightarrow X$ be a holomorphic map that is non-degenerate at some point $x_0\in\C^q$.

  We say that $f$ has {\bf subexponential growth} if the following two conditions are satisfied:

  \vspace{1ex}

  (i)\, there exist constants $C_1>0$ and $r_0>0$ such that \begin{eqnarray}\label{eqn:subexp_1_re}\int\limits_{S_t}|d\tau|_{f^\star\omega}\,d\sigma_{\omega,\,f,\,t}\leq C_1t\, \mbox{Vol}_{\omega,\,f}(B_t), \hspace{5ex} t>r_0;\end{eqnarray}

  (ii)\, for every constant $C>0$, we have: \begin{eqnarray}\label{eqn:subexp_2_re}\limsup\limits_{b\to +\infty}\bigg(\frac{b}{C} - \log F(b)\bigg) = +\infty,\end{eqnarray} where $$F(b):=\int\limits_0^b\mbox{Vol}_{\omega,\,f}(B_t)\,dt = \int\limits_0^b\bigg(\int\limits_{B_t}f^\star\omega_q\bigg)\,dt, \hspace{5ex} b>0.$$

\end{Def}

Definition \ref{Def:p-contact-hyperbolic} introduces a partial hyperbolicity property of $(X,\,\Gamma)$, its ``partial'' feature being reflected in the requirement that $f$ be horizontal. Note that another way of defining horizontality for various maps is well known in the literature and was adopted in [KP23]. Its analogue in our context would consist in requiring that, for every $x\in\C^p$, the image of the differential map $d_xf:\C^p\longrightarrow T_{f(x)}^{1,\,0}X$ of $f$ at $x$ be contained in the fibre ${\cal F}_{\Gamma,\,f(x)}$ at $f(x)$ of ${\cal F}_\Gamma$ whenever $f(x)\in X\setminus\Sigma_\Gamma$ (i.e. whenever ${\cal F}_\Gamma$ is locally free in a neighbourhood of $f(x)$).

Note that, if ${\cal F}_\Gamma$ is supposed locally free (on the whole of $X$), any holomorphic map $f:\C^p\longrightarrow X$ such that $d_xf(\C^p)\subset{\cal F}_{\Gamma,\,f(x)}$ for every $x\in X$ is horizontal in the sense that $f^\star\Gamma = 0$. Indeed, the former condition is equivalent to every $(1,\,0)$-vector field $\xi$ on $X$ that comes from a vector field on $\C^p$ (in the sense that $\xi_{f(x)}$ lies in the image of $d_xf:\C^p\longrightarrow T_{f(x)}^{1,\,0}X$ for every $x\in\C^p$) having the property $\xi\lrcorner\Gamma = 0$. This property means that $\Gamma(\xi,\,\eta_1,\dots , \eta_{p-1}) = 0$ for all $(1,\,0)$-vector fields $\eta_1,\dots , \eta_{p-1}$ on $X$ (not necessarily coming from vector fields on $\C^p$). On the other hand, the condition $f^\star\Gamma = 0$ is equivalent to $\Gamma(\xi_1,\,\xi_2,\dots , \xi_p) = 0$ for all $(1,\,0)$-vector fields $\xi_1,\dots , \xi_p$ on $X$ coming from vector fields on $\C^p$.  

In particular, we conclude that the partial hyperbolicity notion of Definition \ref{Def:p-contact-hyperbolic} is stronger than the one we would obtain if we defined horizontality by means of ${\cal F}_\Gamma$. This ``strength'' is made possible by the hypothesis that a holomorphic $p$-contact structure exist on $X$.

\begin{The}\label{The:partial-p-hyperbolicity_implication} Let $X$ be a compact complex manifold with $\mbox{dim}_\C X = n = 2p+1$. Suppose that $X$ carries a holomorphic $p$-contact structure $\Gamma\in C^\infty_{p,\,0}(X,\,\C)$.

  If $(X,\,\Gamma)$ is {\bf $p$-contact metrically hyperbolic}, then $(X,\,\Gamma)$ is {\bf $p$-contact hyperbolic}.

\end{The}

\noindent {\it Proof.} It is the same as the proof of Theorem 2.5. in [KP23] with the obvious (minor) modifications. For the reader's convenience, we now only point out its main steps. When using Definition \ref{Def:subexp_re} and the notation it employs, we will take the general $q$ to equal $p$.

Suppose $(X,\,\Gamma)$ is $p$-contact metrically hyperbolic and there exists a holomorphic map $f:\C^p\longrightarrow X$ satisfying conditions (i)--(iii) in Definition \ref{Def:p-contact-hyperbolic}. We have to show that this leads to a contradiction.

\vspace{1ex}

{\it Step $1$.} Since $\C^p$ is simply connected, $f$ lifts to $\widetilde{X}$, so there exists a holomorphic map $\widetilde{f}:\C^p\longrightarrow\widetilde{X}$ such that $f=\pi_X\circ\widetilde{f}$. Since $f^\star\Gamma = 0$, we get \begin{eqnarray*}f^\star\omega_p = f^\star\Omega = \widetilde{f}^\star(\pi_X^\star\Omega) = d(\widetilde{f}^\star\widetilde\beta)  \hspace{5ex} \mbox{on}\hspace{1ex} \C^p.\end{eqnarray*}

We notice that $\widetilde{f}^\star\widetilde\beta$ is {\it bounded} in $\C^p$ with respect to the possibly degenerate metric $f^\star\omega$. The proof is the same as that of Lemma 2.6. in [KP23]. (The degeneracy points of $f^\star\omega$ are the points of $\C^p$ where $f^\star\omega$ fails to be positive definite. These are the points of $\C^p$ at which $f$ fails to be an immersion. They form a possibly empty, proper (since $f$ is supposed non-degenerate at least at one point) analytic subset $\Sigma_f$ of $\C^p$.)

\vspace{1ex}

{\it Step $2$.} The classical Fubini Theorem and the H\"older inequality lead to \begin{eqnarray}\label{eqn:Volume_1}\mbox{Vol}_{\omega,\,f}(B_r)\geq 2\,\int\limits_0^r\frac{A_{\omega,\,f}^2(S_t)}{\int\limits_{S_t}|d\tau|_{f^\star\omega}\,d\sigma_{\omega,\,f,\,t}}\,t\,dt, \hspace{5ex} r>0,\end{eqnarray} as in $(9)-(11)$ of [KP23]. This inequality holds for any holomorphic map $f:\C^p\longrightarrow (X,\,\omega)$ that is non-degenerate at some point $x_0\in\C^p$ and takes values in any complex Hermitian manifold. This map need not be either horizontal or of subexponential growth, while the manifold $X$ need not carry a holomorphic $p$-contact structure.

\vspace{1ex}

{\it Step $3$.} An application of the Stokes Theorem and the observations made in {\it Step $1$} lead to \begin{eqnarray}\label{eqn:Volume_Stokes}\nonumber\mbox{Vol}_{\omega,\,f}(B_r) & = & \int\limits_{B_r}f^\star\omega_p = \int\limits_{B_r}f^\star\Omega =\int\limits_{B_r}d(\widetilde{f}^\star\widetilde\beta) = \int\limits_{S_r}\widetilde{f}^\star\widetilde\beta \\
  & \leq & C\,\int\limits_{S_r}d\sigma_{\omega,\,f,\,r} = C\,A_{\omega,\,f}(S_r), \hspace{5ex} r>0,\end{eqnarray} for some constant $C>0$ independent of $r>0$ that comes from the $(f^\star\omega)$-boundedness of $\widetilde{f}^\star\widetilde\beta$. These equalities and inequality are obtained as in (12) of [KP23].

\vspace{1ex}

{\it Step $4$.} The lower, respectively upper, estimate of $\mbox{Vol}_{\omega,\,f}(B_r)$ obtained in {\it Step $2$}, respectively {\it Step $3$} above, together with one part of the subexponential growth assumption on $f$, yield: \begin{eqnarray}\label{eqn:Volume_final}\nonumber\mbox{Vol}_{\omega,\,f}(B_r) & \geq & \frac{2}{C^2}\,\int\limits_0^r \mbox{Vol}_{\omega,\,f}(B_t)\,\frac{t\,\mbox{Vol}_{\omega,\,f}(B_t)}{\int\limits_{S_t}|d\tau|_{f^\star\omega}\,d\sigma_{\omega,\,f,\,t}}\,dt \\
  & \geq & \frac{2}{C_1C^2}\,\int\limits_{r_0}^r\mbox{Vol}_{\omega,\,f}(B_t)\,dt:=\frac{2}{C_1\,C^2}\,F(r), \hspace{5ex} r>r_0\end{eqnarray} for some constants $C_1, r_0>0$, where the last equality constitutes the definition of a function $F:(r_0,\,+\infty)\longrightarrow(0,\,+\infty)$.

  Using this, a Gronwall Lemma-type argument (see end of proof of Theorem 2.5 in [KP23]) shows that $F(a) = 0$ for every $a>r_0$. This means that $\mbox{Vol}_{\omega,\,f}(B_t) = 0$ for every $t>r_0$, which implies $f^\star\omega_p = 0$ on $\C^p$, contradicting the non-degeneracy assumption on $f$ and the property $\omega>0$ on $X$. \hfill $\Box$

\vspace{2ex}

We now use a result from [KP23] on the existence of Ahlfors currents to show that holomorphic $p$-contact manifolds have a natural hyperbolicity property under reasonable extra assumptions. This property is similar to the $p$-contact hyperbolicity introduced in Definition \ref{Def:p-contact-hyperbolic}, but rules out the existence of entire maps from $\C^{p+1}$ (rather than $\C^p$) satisfying a slow growth condition weaker than the {\it subexponential growth} condition. 

\begin{The}\label{The:partial-p+1-hyperbolicity_p-contact} Let $X$ be a compact complex manifold with $\mbox{dim}_\C X = n = 2p+1$ that carries a holomorphic $p$-contact structure $\Gamma\in C^\infty_{p,\,0}(X,\,\C)$.

  Suppose, furthermore, that the sheaves ${\cal F}_\Gamma$ and ${\cal G}_\Gamma$ are both locally free of respective ranks $p+1$ and $p$ and satisfy the condition ${\cal F}_\Gamma\oplus{\cal G}_\Gamma = T^{1,\,0}X$ as vector bundles.

  Then, there exists no holomorphic map $f:\C^{p+1}\longrightarrow X$ with the following three properties:

  \vspace{1ex}

  (i)\, $f$ is non-degenerate at some point $x_0\in\C^{p+1}$;

  \vspace{1ex}

  (ii)\, $f$ satisfies, for some, hence every, Hermitian metric $\omega$, the slow growth condition \begin{eqnarray}\label{eqn:Ahlfors-current}\liminf\limits_{r\to +\infty}\frac{A_{\omega,\,f}(S_r)}{\mbox{Vol}_{\omega,\,f}(B_r)} = 0,\end{eqnarray} where the $(\omega,\,f)$-area $A_{\omega,\,f}(S_r)$ of the Euclidean sphere $S_r\subset\C^{p+1}$ is defined by \begin{eqnarray*}A_{\omega,\,f}(S_r):=\int\limits_{S_r}d\sigma_{\omega,\,f,\,r}>0, \hspace{5ex} r>0,\end{eqnarray*} while the $(\omega,\,f)$-volume of the Euclidean ball $B_r\subset\C^{p+1}$ is defined by \begin{eqnarray*}\mbox{Vol}_{\omega,\,f}(B_r):=\int\limits_{B_r}f^\star\omega_{p+1}.\end{eqnarray*}

  \vspace{1ex}

  (iii)\, $f^\star\Gamma = 0$.

\end{The}

\noindent {\it Proof.} We proceed by contradiction. Suppose a holomorphic map $f:\C^{p+1}\longrightarrow X$ satisfying properties (i), (ii) and (iii) exists.

By Theorem 4.2. in [KP23], any map with the properties (i) and (ii) induces an Ahlfors current $T$ on $X$. This means that $T$ is a strongly semi-positive current of bidimension $(p+1,\,p+1)$ (equivalently, of bidegree $(p,\,p)$) and of mass $1$ with respect to $\omega$ on $X$ having the extra key property that $dT=0$. It is obtained as the weak limit of a subsequence $(T_{r_\nu})_\nu$ of the family of currents \begin{eqnarray*}T_r:=\frac{1}{\mbox{Vol}_{\omega,\,f}(B_r)}\,f_\star[B_r], \hspace{5ex} r>0,\end{eqnarray*} where $[B_r]$ is the current of integration on the ball $B_r\subset\C^{p+1}$.

Using our holomorphic $p$-contact structure $\Gamma$, the $d$-closedness of $T$ and the $\partial\bar\partial$-exactness of the $(p+1,\,p+1)$-form $\partial\Gamma\wedge\bar\partial\overline\Gamma$ imply, via the Stokes theorem, that \begin{eqnarray*}\int\limits_XT\wedge i^{(p+1)^2}\,\partial\Gamma\wedge\bar\partial\overline\Gamma = 0.\end{eqnarray*} Now, since $T$ is a strongly semi-positive current of bidegree $(p,\,p)$ and $i^{(p+1)^2}\,\partial\Gamma\wedge\bar\partial\overline\Gamma$ is a weakly semi-positive $C^\infty$ form of bidegree $(p+1,\,p+1)$, the product $T\wedge i^{(p+1)^2}\,\partial\Gamma\wedge\bar\partial\overline\Gamma$ is a semi-positive current of bidegree $(n,\,n)$ on $X$. Thus, from the vanishing of its integral on $X$ we conclude that \begin{eqnarray}\label{eqn:nn-current_vanishing}T\wedge i^{(p+1)^2}\,\partial\Gamma\wedge\bar\partial\overline\Gamma = 0 \hspace{5ex} \mbox{on}\hspace{1ex} X.\end{eqnarray}

Let $x_0\in\C^{p+1}$ such that $f$ is non-degenerate at $x_0$. Then, $f$ is non-degenerate on a neighbourhood $U$ of $x_0$. Let $\{\xi_1,\dots , \xi_{p+1}\}$ be a holomorphic frame of ${\cal F}_\Gamma$ on a neighbourhood $V$ of $f(x_0)$ in $X$ and let $\{\eta_{p+2},\dots , \eta_n\}$ be a holomorphic frame of ${\cal G}_\Gamma$ on $V$. In particular, $\xi_1\lrcorner\Gamma = \dots = \xi_{p+1}\lrcorner\Gamma = 0$ on $V$.

The property $f^\star\Gamma = 0$ implies that $\tau\lrcorner f^\star\Gamma = 0$ for every vector field $\tau$ on $\C^{p+1}$. In particular, $(f_\star\tau)\lrcorner\Gamma = 0$ at every point of $V$ for every vector field $\tau$ on $\C^{p+1}$. Thus, the tangent vectors $f_\star\tau$ to $X$ at every point in $V$ that come from tangent vectors $\tau$ to $\C^{p+1}$ at points in $U$ are generated by $\xi_1,\dots , \xi_{p+1}$. Since $d_xf:\C^{p+1}\longrightarrow T^{1,\,0}_{f(x)}X$ is injective for every $x\in U$ and the rank of ${\cal F}_\Gamma$ is $p+1$, we deduce that the fibre $({\cal F}_\Gamma)_{f(x)}$, which is already known to be generated by the values at $f(x)$ of $\xi_1,\dots , \xi_{p+1}$, consists precisely of the tangent vectors to $X$ at $f(x)$ that come from tangent vectors to $\C^{p+1}$ at $x$.

On the other hand, the current of integration $[B_r]$ on any ball $B_r\subset\C^{p+1}$ being of bidegree $(0,\,0)$, we have $\tau\lrcorner[B_r] = 0$ and $\overline\tau\lrcorner[B_r] = 0 $ for every $(1,\,0)$-vector field $\tau$ on $\C^{p+1}$, for bidegree reasons.

We deduce that $\xi_j\lrcorner T = 0$ and $\overline\xi_j\lrcorner T = 0$ on $V$ for every $j\in\{1,\dots , p+1\}$. This implies the second equality below, while (\ref{eqn:nn-current_vanishing}) implies the first one: \begin{eqnarray*}0 & = & \overline\xi_{p+1}\lrcorner\dots\lrcorner\overline\xi_1\lrcorner\xi_{p+1}\lrcorner\dots\lrcorner\xi_1\lrcorner\bigg(T\wedge i^{(p+1)^2}\,\partial\Gamma\wedge\bar\partial\overline\Gamma\bigg) \\
  & = & T\wedge\bigg(i^{(p+1)^2}\,(\xi_{p+1}\lrcorner\dots\lrcorner\xi_1\lrcorner\partial\Gamma)\wedge(\overline\xi_{p+1}\lrcorner\dots\lrcorner\overline\xi_1\lrcorner\bar\partial\overline\Gamma)\bigg).\end{eqnarray*} Now, $\xi_{p+1}\lrcorner\dots\lrcorner\xi_1\lrcorner\partial\Gamma$ is a non-vanishing function (which can be made identically $1$ if the local frame $\{\xi_1,\dots , \xi_{p+1}\}$ of ${\cal F}_\Gamma$ is normalised appropriately) on $V$. The same goes for its conjugate function $\overline\xi_{p+1}\lrcorner\dots\lrcorner\overline\xi_1\lrcorner\bar\partial\overline\Gamma$. We infer that $T = 0$ on $V$, hence on $X$, since $x_0\in\C^{p+1}$ was chosen arbitrarily.

The conclusion $T=0$ contradicts the property $\int_X T\wedge\omega_{p+1} = 1$ (which expresses the fact that $T$ has mass $1$ w.r.t. $\omega$ on $X$).  \hfill $\Box$

\section{$p$-contact small deformations}\label{section:contact-def} In this section, we define a type of small deformations of the complex structure of $X$ that preserve a given holomorphic $p$-contact structure $\Gamma$. We then prove unobstructedness to order two.

\vspace{1ex}

\subsection{Meaning of unobstructedness for arbitrary small deformations}\label{subsection:meaning-unobstructedness} We start by recalling a few general, well-known facts in order to fix the notation and the terminology.

Let $X$ be an $n$-dimensional compact complex manifold. A cohomology class $[\theta_1]_{\bar\partial}\in H^{0,\,1}_{\bar\partial}(X,\,T^{1,\,0}X)$ is said to induce a family $(J_t)_{t\in D}$ of small deformations of the complex structure $J_0$ of $X_0:=X$ in the direction of $[\theta_1]_{\bar\partial}$ (where $D$ is a small disc about $0$ in $\C$) if there exists a representative $\psi_1\in C^\infty_{0,\,1}(X,\,T^{1,\,0}X)$ of $[\theta_1]_{\bar\partial}$ and a sequence $(\psi_j)_{j\geq 2}$ of vector-valued forms $\psi_j\in C^\infty_{0,\,1}(X,\,T^{1,\,0}X)$ such that, for all $t\in\C$ sufficiently close to $0$, the power series \begin{eqnarray*}\psi(t)=\psi_1\,t + \psi_2\,t^2 + \dots + \psi_\nu\,t^\nu + \dots\end{eqnarray*} converges (in some topology) and defines a vector-valued form $\psi(t)\in C^\infty_{0,\,1}(X,\,T^{1,\,0}X)$ that satisfies the {\it integrability condition} \begin{eqnarray}\label{eqn:integrability-condition}\bar\partial\psi(t) = \frac{1}{2}\,[\psi(t), \,\psi(t)], \hspace{6ex} t\in D.\end{eqnarray}

  This integrability condition is easily seen to be equivalent to $\bar\partial\psi_1=0$ (always satisfied since $\psi_1$ represents a $\bar\partial$-cohomology class) holding simultaneously with the following sequence of equations: \begin{equation}\label{eqn:integrability-condition_nu}\bar\partial\psi_\nu = \frac{1}{2}\,\sum\limits_{\mu=1}^{\nu-1}[\psi_\mu,\,\psi_{\nu-\mu}] \hspace{3ex} (\mbox{Eq.}\,\,(\nu)), \hspace{3ex} \nu\geq 2. \end{equation}

  The integrability condition (\ref{eqn:integrability-condition}) means that the almost complex structure $J_t$ defined on $X$ by $\psi(t)$ via the operator \begin{eqnarray*}\bar\partial_t\simeq\bar\partial_0 - \psi(t)\end{eqnarray*} is {\it integrable} in the sense that $\bar\partial_t^2 = 0$. (Here, $\bar\partial_0$ is the Cauchy-Riemann operator associated with $J_0$ and the meaning of $\simeq$ is that a $C^\infty$ function $f$ on an open subset of $X$ is holomorphic w.r.t. $J_t$ in the sense that $\bar\partial_t f = 0$ if and only if $(\bar\partial_0 - \psi(t))\,f = 0$.) This amounts to $J_t$ being a complex structure on $X$. In this case, one regards $J_t$ as a deformation of $J_0$.

  Geometrically, the original cohomology class \begin{eqnarray*}[\theta_1]_{\bar\partial} = \bigg[\frac{\partial\psi(t)}{\partial t}_{|t=0}\bigg]_{\bar\partial}\in H^{0,\,1}(X,\,T^{1,\,0}X)\end{eqnarray*} arises, if the integrability condition (\ref{eqn:integrability-condition}) is satisfied, as the tangent vector at $t=0$ to the complex curve $(J_t)_{t\in D}$. Therefore, one regards $[\theta_1]_{\bar\partial}$ as a deformation to order $1$ of $J_0$.

  It is standard to say that the {\it Kuranishi family of $X$ is unobstructed} (to arbitrary order) if every cohomology class $[\theta_1]_{\bar\partial}\in H^{0,\,1}_{\bar\partial}(X,\,T^{1,\,0}X)$ induces a family $(J_t)_{t\in D}$ of small deformations of the complex structure $J_0$ of $X_0:=X$ in the direction of $[\theta_1]_{\bar\partial}$. This is equivalent to the solvability of all the equations in (\ref{eqn:integrability-condition_nu}) for every $[\theta_1]_{\bar\partial}\in H^{0,\,1}_{\bar\partial}(X,\,T^{1,\,0}X)$ and for some representative $\psi_1\in C^\infty_{0,\,1}(X,\,T^{1,\,0}X)$ of the class $[\theta_1]_{\bar\partial}$.

  In this paper, we will also use the following ad hoc definition: the Kuranishi family of $X$ is said to be {\it unobstructed to order two} if (Eq. $(2)$) in (\ref{eqn:integrability-condition_nu}) admits a solution $\psi_2\in C^\infty_{0,\,1}(X,\,T^{1,\,0}X)$ such that $\psi_2\lrcorner u_\Gamma\in\ker\partial$ for every $[\theta_1]_{\bar\partial}\in H^{0,\,1}_{\bar\partial}(X,\,T^{1,\,0}X)$ and for some representative $\psi_1\in C^\infty_{0,\,1}(X,\,T^{1,\,0}X)$ of the class $[\theta_1]_{\bar\partial}$.

  \subsection{Unobstructedness to order two}\label{subsection:order-2-unobstructedness} This is the point of view that we take on deformations in this paper.

  It was observed in [KPU25, Observation 2.1.] that the Fr\"olicher spectral sequence of no holomorphic $p$-contact manifold $X$ degenerates at $E_1$. Meanwhile, the description of the Fr\"olicher spectral sequence of any compact complex manifold $X$ given in [CFGU97] reinterprets the $\C$-vector spaces $E_2^{p,\,q}(X)$ featuring on the second page as the quotients \begin{eqnarray*}E_2^{p,\,q}(X) = \frac{{\cal Z}_2^{p,\,q}(X)}{{\cal C}_2^{p,\,q}(X)}, \hspace{6ex} p,\,q=0,\dots , n,\end{eqnarray*} where the $\C$-vector space ${\cal Z}_2^{p,\,q}(X)$ consists of the $C^\infty$ $(p,\,q)$-forms $u$ such that \begin{eqnarray}\label{eqn:Z_2_def}\bar\partial u=0 \hspace{5ex} \mbox{and} \hspace{5ex} \partial u\in\mbox{Im}\,\bar\partial\end{eqnarray} (these forms $u$ being called {\it $E_2$-closed}), while the $\C$-vector space ${\cal C}_2^{p,\,q}(X)$ consists of the $C^\infty$ $(p,\,q)$-forms $u$ for which there exist forms $\zeta\in C^\infty_{p-1,\,q}(X,\,\C)$ and $\xi\in C^\infty_{p,\,q-1}(X,\,\C)$ such that \begin{eqnarray*}u = \partial\zeta + \bar\partial\xi \hspace{5ex} \mbox{and} \hspace{5ex} \bar\partial\zeta = 0\end{eqnarray*} (these forms $u$ being called {\it $E_2$-exact}).

      On the other hand, the notion of {\it page-$1$-$\partial\bar\partial$-manifold} was introduced in [PSU20a]: a compact complex manifold $X$ is said to have this property if the Fr\"olicher spectral sequence of $X$ degenerates at the second page (a property denoted by $E_2(X) = E_\infty(X)$) and the De Rham cohomology of $X$ is {\it pure}. In this paper, we will use the following consequence of Theorem 4.3. in [PSU20b]: if $X$ is a {\it page-$1$-$\partial\bar\partial$-manifold}, the following equality of $\C$-vector spaces holds in every bidegree $(p,\,q)$: \begin{eqnarray}\label{eqn:ddbar_im_Z-2}\partial\bigg({\cal Z}_2^{p,\,q}(X)\bigg) = \mbox{Im}\,(\partial\bar\partial).\end{eqnarray} The inclusion ``$\supset$'' holds on any $X$, it is the inclusion ``$\subset$'' that follows from the {\it page-$1$-$\partial\bar\partial$}-assumption.

\vspace{2ex}      

 The class of small deformations of the complex structure of a $p$-contact manifold that we introduce and study in this paper is described in the following

\begin{Def}\label{Def:contact_01_T10} Let $X$ be a compact complex manifold with $\mbox{dim}_\C X = n = 2p+1$. Suppose there exists a holomorphic $p$-contact structure $\Gamma\in C^\infty_{p,\,0}(X,\,\C)$ on $X$. 

(i)\, The $\C$-vector subspace \begin{eqnarray}\label{eqn:contact_01_T10}\nonumber H^{0,\,1}_{[\Gamma],\,def}(X,\,T^{1,\,0}X) & := & \bigg\{[\theta]_{\bar\partial}\in H^{0,\,1}_{\bar\partial}(X,\,T^{1,\,0}X)\,\mid\,[\theta\lrcorner\Gamma]_{\bar\partial} = 0\in H^{p-1,\,1}_{\bar\partial}(X,\,\C) \hspace{1ex} \mbox{and} \hspace{1ex} \theta\lrcorner\partial\Gamma\in{\cal Z}_2^{p,\,1}(X)\bigg\} \\
    & \subset &  H^{0,\,1}_{\bar\partial}(X,\,T^{1,\,0}X)\end{eqnarray} is called the space of (infinitesimal) {\bf $p$-contact deformations} of $X$.

\vspace{1ex}

(ii)\, We say that the $p$-contact deformations of $X$ are {\bf unobstructed to order two} if every class $[\theta_1]_{\bar\partial}\in H^{0,\,1}_{[\Gamma],\,def}(X,\,T^{1,\,0}X)$ has a representative $\psi_1\in C^\infty_{0,\,1}(X,\,T^{1,\,0}X)$ such that (Eq. $(2)$) in (\ref{eqn:integrability-condition_nu}) admits a solution $\psi_2\in C^\infty_{0,\,1}(X,\,T^{1,\,0}X)$ with the property $\psi_2\lrcorner u_\Gamma\in\ker\partial$.

\end{Def}

\begin{Rem} Definition \ref{Def:contact_01_T10} is correct in the sense that the space $H^{0,\,1}_{[\Gamma],\,def}(X,\,T^{1,\,0}X)$ is {\bf well defined}.

\end{Rem}

\noindent {\it Proof.} Since $\bar\partial(\theta\lrcorner\Gamma) = (\bar\partial\theta)\lrcorner\Gamma + \theta\lrcorner(\bar\partial\Gamma)$ (see (\ref{eqn:theta_contraction-Leibniz})) and since $\bar\partial\theta = 0$ (for any $\theta\in C^\infty_{0,\,1}(X,\,T^{1,\,0}X)$ representing a $\bar\partial$-cohomology class) and $\bar\partial\Gamma=0$ (by hypothesis), we get $\bar\partial(\theta\lrcorner\Gamma) = 0$, so the $\bar\partial$-cohomology class $[\theta\lrcorner\Gamma]_{\bar\partial}$ is well defined as an element of $H^{p-1,\,1}_{\bar\partial}(X,\,\C)$. Moreover, from (\ref{eqn:xi_contraction-Leibniz}) we get $\bar\partial(\xi\lrcorner\Gamma) = (\bar\partial\xi)\lrcorner\Gamma - \xi\lrcorner(\bar\partial\Gamma) = (\bar\partial\xi)\lrcorner\Gamma$, so the class $[\theta\lrcorner\Gamma]_{\bar\partial}$ is independent of the choice of representative of the class $[\theta]_{\bar\partial}$. Thus, the condition $[\theta\lrcorner\Gamma]_{\bar\partial} = 0$ in (\ref{eqn:contact_01_T10}) is well posed as it depends only on $[\theta]_{\bar\partial}$.

Similarly, the condition $\theta\lrcorner\partial\Gamma\in{\cal Z}_2^{p,\,1}(X)$ is well posed. Indeed, if $\theta$ is changed to another representative $\theta + \bar\partial\xi$ (with $\xi\in C^\infty(X,\,T^{1,\,0}X$)) of its $\bar\partial$-cohomology class, we get: $(\theta + \bar\partial\xi)\lrcorner\partial\Gamma = \theta\lrcorner\partial\Gamma + \bar\partial(\xi\lrcorner\partial\Gamma)$ (because $\bar\partial\partial\Gamma = 0$). Since any $\bar\partial$-exact form lies in ${\cal Z}_2^{p,\,1}(X)$ (as can be checked right away from (\ref{eqn:Z_2_def})), we infer that $(\theta + \bar\partial\xi)\lrcorner\partial\Gamma\in{\cal Z}_2^{p,\,1}(X)$ if and only if $\theta\lrcorner\partial\Gamma\in{\cal Z}_2^{p,\,1}(X)$. \hfill $\Box$

\vspace{1ex}


\begin{Lem}\label{Lem:contact_2p_1_scalar} The image of $H^{0,\,1}_{[\Gamma],\,def}(X,\,T^{1,\,0}X)\subset  H^{0,\,1}_{\bar\partial}(X,\,T^{1,\,0}X)$ under the Calabi-Yau isomorphism $T_{[\Gamma]}$ of (\ref{eqn:C_Y-isomorphism_cohom}) is the vector subspace \begin{eqnarray}\label{eqn:contact_2p_1_scalar}H^{2p,\,1}_{[\Gamma],\,def}(X,\,\C):=\bigg\{[\Gamma\wedge(\theta\lrcorner\partial\Gamma)]_{\bar\partial}\in H^{2p,\,1}_{\bar\partial}(X,\,\C)\,\mid\,[\theta]_{\bar\partial}\in H^{0,\,1}_{[\Gamma],\,def}(X,\,T^{1,\,0}X)\bigg\}\subset H^{2p,\,1}_{\bar\partial}(X,\,\C).\end{eqnarray}

\end{Lem}

\noindent {\it Proof.} Indeed, $n-1 = 2p$, while the equality \begin{eqnarray*}\theta\lrcorner(\Gamma\wedge\partial\Gamma) = (\theta\lrcorner\Gamma)\wedge\partial\Gamma + \Gamma\wedge(\theta\lrcorner\partial\Gamma),\end{eqnarray*} the property $[\theta\lrcorner\Gamma]_{\bar\partial} = 0$ for every $[\theta]_{\bar\partial}\in H^{0,\,1}_{[\Gamma]}(X,\,T^{1,\,0}X)$ and $\bar\partial(\partial\Gamma) = 0$ imply that $(\theta\lrcorner\Gamma)\wedge\partial\Gamma\in\mbox{Im}\,\bar\partial$ and that $[\theta\lrcorner(\Gamma\wedge\partial\Gamma)]_{\bar\partial} = [\Gamma\wedge(\theta\lrcorner\partial\Gamma)]_{\bar\partial}$.  \hfill $\Box$

\vspace{2ex}

Our partial unobstructedness result is the following

\begin{The}\label{The:partial-unobstructedness} Let $X$ be a compact complex manifold with $\mbox{dim}_\C X = n = 2p+1$. Suppose that $X$ is a {\bf page-$1$-$\partial\bar\partial$-manifold} that carries a {\bf holomorphic $p$-contact structure} $\Gamma\in C^\infty_{p,\,0}(X,\,\C)$ such that the sheaf ${\cal F}_\Gamma$ is {\bf cohomologically integrable in bidegree $(0,\,1)$}.

  Then, the {\bf $p$-contact deformations} of $X$ are {\bf unobstructed to order two}.

\end{The}

\noindent {\it Proof.} Let $[\theta_1]_{\bar\partial}\in H^{0,\,1}_{[\Gamma],\,def}(X,\,T^{1,\,0}X)$ be arbitrary. We will prove that $[\theta_1]_{\bar\partial}$ induces a family of deformations of the complex structure of $X$ in the direction of $[\theta_1]_{\bar\partial}$.

\vspace{1ex}

$\bullet$ {\it Step 1.} Let $\theta_1\in C^\infty_{0,\,1}(X,\,T^{1,\,0}X)$ be an arbitrary representative of the class $[\theta_1]_{\bar\partial}$. Then, $\bar\partial\theta_1 = 0$ and, as noticed in Lemma \ref{Lem:contact_2p_1_scalar}, we have \begin{eqnarray*}T_{[\Gamma]}([\theta_1]_{\bar\partial}) = [\Gamma\wedge(\theta_1\lrcorner\partial\Gamma)]_{\bar\partial}\in H^{2p,\,1}_{\bar\partial}(X,\,\C).\end{eqnarray*}

Let $\eta_1:= \Gamma\wedge(\theta_1\lrcorner\partial\Gamma)\in C^\infty_{2p,\,1}(X,\,\C)$. By the definition (\ref{eqn:contact_01_T10}) of $H^{0,\,1}_{[\Gamma]}(X,\,T^{1,\,0}X)$, we have $\theta_1\lrcorner\partial\Gamma\in{\cal Z}_2^{p,\,1}(X)$. This means (cf. (\ref{eqn:Z_2_def})) that $\bar\partial(\theta_1\lrcorner\partial\Gamma) = 0$ and $\partial(\theta_1\lrcorner\partial\Gamma)\in\mbox{Im}\,\bar\partial$. We will now observe that these two properties, together with $\bar\partial\Gamma = 0$, imply \begin{eqnarray}\label{eqn:eta_1_Z-2}\eta_1\in{\cal Z}_2^{2p,\,1}(X).\end{eqnarray}

Indeed, $\bar\partial\eta_1 = \bar\partial\Gamma\wedge(\theta_1\lrcorner\partial\Gamma) + (-1)^p\,\Gamma\wedge\bar\partial(\theta_1\lrcorner\partial\Gamma) = 0$. Meanwhile, $\partial\Gamma\wedge\partial\Gamma = 0$ for bidegree reasons (since this is a form of bidegree $(2p+2,\,0) = (n+1,\,0$)). Then, \begin{eqnarray*}0 = \theta_1\lrcorner(\partial\Gamma\wedge\partial\Gamma) = (\theta_1\lrcorner\partial\Gamma)\wedge\partial\Gamma + \partial\Gamma\wedge(\theta_1\lrcorner\partial\Gamma) = 2\,\partial\Gamma\wedge(\theta_1\lrcorner\partial\Gamma),\end{eqnarray*} where the last equality follows from $p+1$ being even (since $p$ is odd). We get $\partial\Gamma\wedge(\theta_1\lrcorner\partial\Gamma) = 0$ and \begin{eqnarray*}\partial\eta_1 = \partial\Gamma\wedge(\theta_1\lrcorner\partial\Gamma) + (-1)^p\,\Gamma\wedge\partial(\theta_1\lrcorner\partial\Gamma) = (-1)^p\,\Gamma\wedge\partial(\theta_1\lrcorner\partial\Gamma).\end{eqnarray*} Since $\Gamma$ is $\bar\partial$-closed and $\partial(\theta_1\lrcorner\partial\Gamma)$ is $\bar\partial$-exact, $\Gamma\wedge\partial(\theta_1\lrcorner\partial\Gamma)$ is $\bar\partial$-exact. This proves that $\partial\eta_1$ is $\bar\partial$-exact.

Thus, (\ref{eqn:eta_1_Z-2}) is proven. We infer that \begin{eqnarray*}\partial\eta_1\in\partial\bigg({\cal Z}_2^{2p,\,1}(X)\bigg) = \mbox{Im}\,(\partial\bar\partial).\end{eqnarray*} (See (\ref{eqn:ddbar_im_Z-2}) for the last equality, a consequence of the {\it page-$1$-$\partial\bar\partial$}-assumption on $X$.) Consequently, there exists $\zeta_1\in C^\infty_{n-1,\,0}(X,\,T^{1,\,0}X)$ such that $\partial\eta_1 = \partial\bar\partial\zeta_1$. Thanks to the Calabi-Yau isomorphism $T_\Gamma$, there exists a unique $\psi_1\in C^\infty_{0,\,1}(X,\,T^{1,\,0}X)$ such that \begin{eqnarray}\label{eqn:psi_1_def}\psi_1\lrcorner u_\Gamma = \eta_1 - \bar\partial\zeta_1.\end{eqnarray}

By construction, $\bar\partial(\psi_1\lrcorner u_\Gamma) = 0$ (hence also $\bar\partial\psi_1 = 0$) and $\partial(\psi_1\lrcorner u_\Gamma) = 0$. Moreover, \begin{eqnarray*}[\psi_1\lrcorner u_\Gamma]_{\bar\partial} = [\eta_1]_{\bar\partial} = [\theta_1\lrcorner u_\Gamma]_{\bar\partial}.\end{eqnarray*} In particular, $[\psi_1]_{\bar\partial} = [\theta_1]_{\bar\partial}$. It is in order to achieve the $\partial$-closedness (in addition to the $\bar\partial$-closedness) of $\psi_1\lrcorner u_\Gamma$ that we changed the representative $\theta_1$ to $\psi_1$ of the original class $[\theta_1]_{\bar\partial}\in H^{0,\,1}_{[\Gamma]}(X,\,T^{1,\,0}X)$. This $\partial$-closedness will enable us to apply the Tian-Todorov Lemma \ref{Lem:Tian-Todorov} at the next stage. 

\vspace{1ex}

$\bullet$ {\it Step 2.} Now that $\psi_1$ has been constructed, we need to solve the equation \begin{equation}\label{eqn:integrability-condition_2}\bar\partial\psi_2 = \frac{1}{2}\,[\psi_1,\,\psi_1] \hspace{3ex} (\mbox{Eq.}\,\,(2)) \end{equation} by finding a solution $\psi_2\in C^\infty_{0,\,1}(X,\,T^{1,\,0}X)$ such that $\psi_2\lrcorner u_\Gamma\in C^\infty_{n-1,\,1}(X,\,\C)$ is $\partial$-closed. We will actually achieve rather more by proving the existence of a solution $\psi_2$ such that $\psi_2\lrcorner u_\Gamma$ is $\partial$-exact.

Since $\partial(\psi_1\lrcorner u_\Gamma) = 0$, the Tian-Todorov Lemma \ref{Lem:Tian-Todorov} yields \begin{equation}\label{eqn:T-T_stage-2}[\psi_1,\,\psi_1]\lrcorner u_\Gamma = -\partial\bigg(\psi_1\lrcorner(\psi_1\lrcorner u_\Gamma)\bigg).\end{equation}

We will now prove that $[\psi_1,\,\psi_1]\lrcorner u_\Gamma$ is $\bar\partial$-exact. Since $u_\Gamma = \Gamma\wedge\partial\Gamma$, we have \begin{equation}\label{eqn:psi-1_psi-1_bracket_u-Gamma}[\psi_1,\,\psi_1]\lrcorner u_\Gamma = \bigg([\psi_1,\,\psi_1]\lrcorner\Gamma\bigg)\wedge\partial\Gamma + (-1)^p\,\Gamma\wedge\bigg([\psi_1,\,\psi_1]\lrcorner\partial\Gamma\bigg).\end{equation}

Now, recall that $[\theta_1\lrcorner\Gamma]_{\bar\partial} = [\psi_1\lrcorner\Gamma]_{\bar\partial} = 0\in H^{p-1,\,1}_{\bar\partial}(X,\,\C)$ since $[\psi_1]_{\bar\partial} = [\theta_1]_{\bar\partial} \in H^{0,\,1}_{[\Gamma]}(X,\,T^{1,\,0}X)$. Hence, by the cohomological integrability in bidegree $(0,\,1)$ assumption on ${\cal F}_\Gamma$, there exists $\beta\in C^\infty_{p-1,\,1}(X,\,\C)$ such that \begin{eqnarray}\label{eqn:psi-1_psi-1_Gamma_potential}[\psi_1,\,\psi_1]\lrcorner\Gamma = \bar\partial\beta.\end{eqnarray} Since $\partial\Gamma\in\ker\bar\partial$, this implies that \begin{equation}\label{eqn:psi-1_psi-1_bracket_u-Gamma_1}\bigg([\psi_1,\,\psi_1]\lrcorner\Gamma\bigg)\wedge\partial\Gamma = \bar\partial(\beta\wedge\partial\Gamma).\end{equation}

It remains to show that the last term in (\ref{eqn:psi-1_psi-1_bracket_u-Gamma}) is $\bar\partial$-exact. Since $\psi_1\lrcorner\Gamma$ is $\bar\partial$-exact, there exists $\alpha\in C^\infty_{p-1,\,0}(X,\,\C)$ such that \begin{eqnarray}\label{eqn:psi-1_Gamma_potential}\psi_1\lrcorner\Gamma = \bar\partial\alpha.\end{eqnarray}

\begin{Claim}\label{Claim:psi-1_psi-1_Gamma} The following identity holds: \begin{eqnarray*}\label{eqn:psi-1_psi-1_Gamma}[\psi_1,\,\psi_1]\lrcorner\Gamma + \psi_1\lrcorner(\psi_1\lrcorner\partial\Gamma) = \bar\partial\bigg(\partial(\psi_1\lrcorner\alpha) - 2\psi_1\lrcorner\partial\alpha\bigg).\end{eqnarray*}

\end{Claim}

\noindent {\it Proof of Claim.} The general identity (\ref{eqn:Lie-derivatives-mixed_prop_2}) yields the first equality below: \begin{eqnarray*}[\psi_1,\,\psi_1]\lrcorner\Gamma + \psi_1\lrcorner(\psi_1\lrcorner\partial\Gamma) & = & -[L_{\psi_1},\,\psi_1\lrcorner\cdot\,](\Gamma) + \psi_1\lrcorner(\psi_1\lrcorner\partial\Gamma) = -L_{\psi_1}(\psi_1\lrcorner\Gamma) + \psi_1\lrcorner L_{\psi_1}\Gamma + \psi_1\lrcorner(\psi_1\lrcorner\partial\Gamma) \\
  & = & -L_{\psi_1}(\bar\partial\alpha) + \psi_1\lrcorner L_{\psi_1}\Gamma + \psi_1\lrcorner(\psi_1\lrcorner\partial\Gamma),\end{eqnarray*} where the last equality follows from (\ref{eqn:psi-1_Gamma_potential}). Now, $\partial\bar\partial\alpha = \partial(\psi_1\lrcorner\Gamma) = L_{\psi_1}\Gamma + \psi_1\lrcorner\partial\Gamma$ (see (\ref{eqn:Lie-derivatives-mixed_def}) for the last equality). Applying $\psi_1\lrcorner\cdot$ on either side of this equality, we get  \begin{eqnarray*}\psi_1\lrcorner L_{\psi_1}\Gamma + \psi_1\lrcorner(\psi_1\lrcorner\partial\Gamma) = \psi_1\lrcorner\partial\bar\partial\alpha.\end{eqnarray*} Putting together these pieces of information, we get: \begin{eqnarray*}[\psi_1,\,\psi_1]\lrcorner\Gamma + \psi_1\lrcorner(\psi_1\lrcorner\partial\Gamma) & = & -L_{\psi_1}(\bar\partial\alpha) + \psi_1\lrcorner\partial\bar\partial\alpha = -L_{\psi_1}(\bar\partial\alpha) - \psi_1\lrcorner\bar\partial(\partial\alpha) \\
  & = & -L_{\psi_1}(\bar\partial\alpha) - \bar\partial\bigg(\psi_1\lrcorner\partial\alpha\bigg),\end{eqnarray*} where we used the property $\bar\partial\psi_1 = 0$ to get the last equality.

Meanwhile, thanks to (\ref{eqn:Lie-derivatives-mixed_prop_0}) and to $\bar\partial\psi_1 = 0$, we have $-L_{\psi_1}(\bar\partial\alpha) = \bar\partial(L_{\psi_1}\alpha)$. Thus, the above equality translates to \begin{eqnarray}\label{eqn:psi-1_psi-1_Gamma_proof_2}\nonumber[\psi_1,\,\psi_1]\lrcorner\Gamma + \psi_1\lrcorner(\psi_1\lrcorner\partial\Gamma) = \bar\partial\bigg(L_{\psi_1}\alpha - \psi_1\lrcorner\partial\alpha\bigg) = \bar\partial\bigg(\partial(\psi_1\lrcorner\alpha) - 2\psi_1\lrcorner\partial\alpha\bigg),\end{eqnarray} where the last equality follows from $\partial(\psi_1\lrcorner\alpha) = L_{\psi_1}\alpha + \psi_1\lrcorner\partial\alpha$ (see (\ref{eqn:Lie-derivatives-mixed_def})).

This proves Claim \ref{Claim:psi-1_psi-1_Gamma}.  \hfill $\Box$

\begin{Claim}\label{Claim:psi-1_psi-1_del-Gamma} The following identity holds: \begin{eqnarray*}[\psi_1,\,\psi_1]\lrcorner\partial\Gamma = 2\, \psi_1\lrcorner\partial(\psi_1\lrcorner\partial\Gamma) - \partial\bigg(\psi_1\lrcorner(\psi_1\lrcorner\partial\Gamma)\bigg).\end{eqnarray*}

\end{Claim}

\noindent {\it Proof of Claim.} The general identity (\ref{eqn:Lie-derivatives-mixed_prop_2}) yields the first equality below: \begin{eqnarray*}[\psi_1,\,\psi_1]\lrcorner\partial\Gamma & = & -[L_{\psi_1},\,\psi_1\lrcorner\cdot\,](\partial\Gamma) = -L_{\psi_1}(\psi_1\lrcorner\partial\Gamma) + \psi_1\lrcorner L_{\psi_1}(\partial\Gamma).\end{eqnarray*}

Now, $\partial(\psi_1\lrcorner\partial\Gamma) = L_{\psi_1}(\partial\Gamma) + \psi_1\lrcorner\partial^2\Gamma = L_{\psi_1}(\partial\Gamma)$. Applying $\psi_1\lrcorner\cdot$ on either side, we get  \begin{eqnarray*}\psi_1\lrcorner L_{\psi_1}(\partial\Gamma) = \psi_1\lrcorner\partial(\psi_1\lrcorner\partial\Gamma).\end{eqnarray*}

Meanwhile, \begin{eqnarray*}\partial\bigg(\psi_1\lrcorner(\psi_1\lrcorner\partial\Gamma)\bigg) = L_{\psi_1}(\psi_1\lrcorner\partial\Gamma) + \psi_1\lrcorner\partial(\psi_1\lrcorner\partial\Gamma).\end{eqnarray*}

Plugging the expressions for $\psi_1\lrcorner L_{\psi_1}(\partial\Gamma)$ and $-L_{\psi_1}(\psi_1\lrcorner\partial\Gamma)$ given by the last two identities into the expression for $[\psi_1,\,\psi_1]\lrcorner\partial\Gamma$, we get the claimed identity.    \hfill $\Box$

\vspace{2ex}

We now plug the expression for $\psi_1\lrcorner(\psi_1\lrcorner\partial\Gamma)$ given by Claim \ref{Claim:psi-1_psi-1_Gamma} into the formula proved in Claim \ref{Claim:psi-1_psi-1_del-Gamma}. We get: \begin{eqnarray}\label{eqn:psi-1_psi-1_del-Gamma}\nonumber[\psi_1,\,\psi_1]\lrcorner\partial\Gamma & = & 2\, \psi_1\lrcorner\partial(\psi_1\lrcorner\partial\Gamma) + \partial\bigg([\psi_1,\,\psi_1]\lrcorner\Gamma\bigg) + \bar\partial\partial\bigg(\partial(\psi_1\lrcorner\alpha) - 2\psi_1\lrcorner\partial\alpha\bigg) \\
  & = & 2\, \psi_1\lrcorner\partial(\psi_1\lrcorner\partial\Gamma) - \bar\partial\partial\beta - 2\bar\partial\partial\bigg(\psi_1\lrcorner\partial\alpha\bigg),\end{eqnarray} where for the last equality we used (\ref{eqn:psi-1_psi-1_Gamma_potential}).

Now, recall that $[\psi_1]_{\bar\partial} = [\theta_1]_{\bar\partial}\in H^{0,\,1}_{[\Gamma]}(X,\,T^{1,\,0}X)$. Hence, by definition (\ref{eqn:contact_01_T10}) of $H^{0,\,1}_{[\Gamma]}(X,\,T^{1,\,0}X)$, we have $\psi_1\lrcorner\partial\Gamma\in{\cal Z}_2^{p,\,1}(X)$. Thanks to (\ref{eqn:Z_2_def}), this implies that $\partial(\psi_1\lrcorner\partial\Gamma)\in\mbox{Im}\,\bar\partial$. Consequently, there exists $\delta\in C^\infty_{p+1,\,0}(X,\,\C)$ such that \begin{eqnarray}\label{eqn:del_psi-1_del-Gamma_potential}\partial(\psi_1\lrcorner\partial\Gamma) = \bar\partial\delta.\end{eqnarray}

Hence, $\psi_1\lrcorner\partial(\psi_1\lrcorner\partial\Gamma) = \psi_1\lrcorner\bar\partial\delta = \bar\partial(\psi_1\lrcorner\delta)$. Therefore, (\ref{eqn:psi-1_psi-1_del-Gamma}) becomes \begin{eqnarray}\label{eqn:psi-1_psi-1_del-Gamma_final}[\psi_1,\,\psi_1]\lrcorner\partial\Gamma = \bar\partial\bigg(2\,\psi_1\lrcorner\delta - \partial\beta - 2\,\partial(\psi_1\lrcorner\partial\alpha)\bigg).\end{eqnarray} Since $\bar\partial\Gamma = 0$, this shows in particular that the last term in (\ref{eqn:psi-1_psi-1_bracket_u-Gamma}) is $\bar\partial$-exact.

Putting together (\ref{eqn:psi-1_psi-1_bracket_u-Gamma_1}) and (\ref{eqn:psi-1_psi-1_del-Gamma_final}), (\ref{eqn:psi-1_psi-1_bracket_u-Gamma}) becomes: \begin{equation}\label{eqn:psi-1_psi-1_bracket_u-Gamma_final}[\psi_1,\,\psi_1]\lrcorner u_\Gamma = \bar\partial\bigg(\beta\wedge\partial\Gamma + \Gamma\wedge\bigg(2\,\psi_1\lrcorner\delta - \partial\beta - 2\,\partial(\psi_1\lrcorner\partial\alpha)\bigg)\bigg)\in\mbox{Im}\,\bar\partial.\end{equation}

Combined with (\ref{eqn:T-T_stage-2}) and with the obvious fact that $\psi_1\lrcorner(\psi_1\lrcorner u_\Gamma)\in\ker\bar\partial$, this proves that $\psi_1\lrcorner(\psi_1\lrcorner u_\Gamma)\in{\cal Z}_2^{n-2,\,2}(X)$. In particular, (\ref{eqn:T-T_stage-2}) reads: \begin{equation*}[\psi_1,\,\psi_1]\lrcorner u_\Gamma = -\partial\bigg(\psi_1\lrcorner(\psi_1\lrcorner u_\Gamma)\bigg)\in\partial\bigg({\cal Z}_2^{n-2,\,2}(X)\bigg) = \mbox{Im}\,(\partial\bar\partial).\end{equation*} (See (\ref{eqn:ddbar_im_Z-2}) for the last equality, a consequence of the {\it page-$1$-$\partial\bar\partial$}-assumption on $X$.) Therefore, there exists $\Phi_2\in C^\infty_{n-2,\,1}(X,\,\C)$ such that \begin{equation}\label{eqn:eqn-2_dd-bar_solution}\bar\partial\partial\Phi_2 = \frac{1}{2}\,[\psi_1,\,\psi_1]\lrcorner u_\Gamma.\end{equation}

Now, $\partial\Phi_2\in C^\infty_{n-1,\,1}(X,\,\C)$, so the Calabi-Yau isomorphism $T_\Gamma$ ensures the existence of a unique $\psi_2\in C^\infty_{0,\,1}(X,\,T^{1,\,0}X)$ such that $\psi_2\lrcorner u_\Gamma = \partial\Phi_2$. In particular, (\ref{eqn:eqn-2_dd-bar_solution}) translates to \begin{equation}\label{eqn:eqn-2_d-bar_solution}\bar\partial(\psi_2\lrcorner u_\Gamma) = \frac{1}{2}\,[\psi_1,\,\psi_1]\lrcorner u_\Gamma.\end{equation} This means that $\psi_2$ is a solution to equation (\ref{eqn:integrability-condition_2}) (also called $(\mbox{Eq.}\,\,(2))$) with the property that $\psi_2\lrcorner u_\Gamma\in\mbox{Im}\,\partial\subset\ker\partial$. \hfill $\Box$

\vspace{3ex}

\noindent {\bf Acknowledgements.} The second-named author visited his co-authors separately in their respective institutions and is grateful to the Osaka University and to the University of Zaragoza for support and hospitality.  The third-named author  was partially supported by grant PID2023-148446NB-I00, funded by MICIU/AEI/10.13039/501100011033, 
and by grant E22-23R ``Algebra y Geometr\' ia'' (Gobierno de Arag\'on/FEDER).

\vspace{3ex}

\noindent {\bf References.} \\






\vspace{1ex}

\noindent [Cle05]\, H. Clemens --- {\it Geometry of Formal Kuranishi Theory} --- Adv. Math., {\bf 198} (2005) 311-365.

\vspace{1ex}

\noindent [CFGU97]\, L.A. Cordero, M. Fern\'andez, A.Gray, L. Ugarte --- {\it A General Description of the Terms in the Fr\"olicher Spectral Sequence} --- Differential Geom. Appl. {\bf 7} (1997), no. 1, 75--84.

\vspace{1ex}

\noindent [Dem02]\, J.-P. Demailly --- {\it On the Frobenius Integrability of Certain Holomorphic $p$-Forms } --- In: Bauer, I., Catanese, F., Peternell, T., Kawamata, Y., Siu, YT. (eds) Complex Geometry. Springer, Berlin, Heidelberg. \url{https://doi.org/10.1007/978-3-642-56202-0_6}

\vspace{1ex}

\noindent [Dem97]\, J.-P. Demailly --- {\it Complex Analytic and Algebraic Geometry} --- \url{https://www-fourier.univ-grenoble-alpes.fr/~demailly/manuscripts/agbook.pdf}



\vspace{1ex}

\noindent [FOU15]\, A. Fino, A. Otal, L. Ugarte --- {\it Six-dimensional solvmanifolds with holomorphically trivial canonical bundle} --- Int. Math. Res. Not. IMRN(2015), no. 24, 13757--13799, doi: 10.1093/imrn/rnv112

\vspace{1ex}

\noindent [FRR19]\, A. Fino, S. Rollenske, J. Ruppenthal --- {\it Dolbeault cohomology of complex nilmanifolds foliated in toroidal groups} --- Q. J. Math. {\bf 70} (2019), no.4, 1265--1279.

\vspace{1ex}

\noindent [HK74]\, R. Harvey, A.W. Knapp -- {\it Positive $(p,\,p)$-forms, Wirtinger's Inequality and Currents} -- In {\it Value distribution theory (Proc. Tulane Univ. Program,Tulane Univ., New Orleans, La., 1972-1973), Part A}, p. 43-62, Dekker, New York, 1974.


\vspace{1ex}

\noindent [KP23]\, H. Kasuya, D. Popovici --- {\it Partially Hyperbolic Compact Complex Manifolds} --- 
Rev. Mat. Iberoam. {\bf 41} (2025), no. 5, pp. 1795--1832; 
doi 10.4171/RMI/1555

\vspace{1ex}

\noindent [KPU25]\, H. Kasuya, D. Popovici, L. Ugarte --- {\it Higher-Degree Holomorphic Contact Structures} --- 
arXiv:2502.01447v2



\vspace{1ex}

\noindent [Kod86]\, K. Kodaira --- {\it Complex Manifolds and Deformations of Complex Structures} --- Grundlehren der Math. Wiss. {\bf 283}, Springer (1986).



\vspace{1ex}

\noindent [Ma24]\, Y. Ma --- {\it Strongly Gauduchon Hyperbolicity and Two Other Types of Hyperbolicity} --- arXiv e-print DG 2404.08830.

\vspace{1ex}

\noindent [Mar23a]\, S. Marouani --- {\it SKT Hyperbolic and Gauduchon Hyperbolic Compact Complex Manifolds} --- 
arXiv e-print DG 2305.08122; 
Bull. Soc. Math. France {\bf 152} (3) (2024), 519--545; 
doi: 10.24033/bsmf.2894

\vspace{1ex}

\noindent [Mar23b]\, S. Marouani --- {\it A Comparative Study between K\"ahler and Non-K\"ahler Hyperbolicity} --- arXiv e-print CV 2306.15036; 
Complex Var. Elliptic Equ. (2025), 1--10; doi: 10.1080/17476933. 2025.2532508

\vspace{1ex}

\noindent [MP22a]\, S. Marouani, D. Popovici --- {\it Balanced Hyperbolic and Divisorially Hyperbolic Compact Complex Manifolds} --- arXiv e-print CV 2107.08972v2;  Math. Res. Lett., Vol. 30, No. 6, 1813--1855 (2023).

\vspace{1ex}

\noindent [MP22b]\, S. Marouani, D. Popovici --- {\it Some Properties of Balanced Hyperbolic Compact Complex Manifolds } --- arXiv e-print CV 2107.09522;  Internat. J. Math., Vol. 33, No. 3 (2022) 2250019.

\vspace{1ex}

\noindent [Nak75]\, I. Nakamura --- {\it Complex parallelisable manifolds and their small deformations} --- J. Differ. Geom. {\bf 10}, (1975), 85--112.

\vspace{1ex}

\noindent [OU23-25]\, A. Otal, L. Ugarte --- {\it Six dimensional homogeneous spaces with holomorphically trivial canonical bundle} --- J. Geom. Phys. {\bf 194} (2023), 105014, 27 pp. 
[Addendum: J. Geom. Phys. {\bf 216} (2025), 105570, 5 pp.].



\vspace{1ex}

\noindent [Pop18]\, D. Popovici --- {\it Non-K\"ahler Mirror Symmetry of the Iwasawa Manifold} --- Int. Math. Res. Not. IMRN 2020, no. {\bf 23}, 9471--9538, doi:10.1093/imrn/rny256.

\vspace{1ex}

\noindent [Pop19]\,  D. Popovici --- {\it Holomorphic Deformations of Balanced Calabi-Yau $\partial\bar\partial$-Manifolds} --- Ann. Inst. Fourier, {\bf 69} (2019) no. 2, pp. 673--728. doi : 10.5802/aif.3254.

\vspace{1ex}

\noindent [Pop25]\, D. Popovici --- {\it Generalised Hermite-Einstein Fibre Metrics and Slope Stability for Holomorphic Vector Bundles} --- arXiv:2512.24932v2 [math.AG].

\vspace{1ex}

\noindent [PSU20a]\,  D. Popovici, J. Stelzig, L. Ugarte --- {\it Higher-Page Hodge Theory of Compact Complex Manifolds} --- arXiv e-print AG 2001.02313v3; Ann. Sc. Norm. Super. Pisa Cl. Sci. (5) Vol. XXV (2024), 1431--1464. 

\vspace{1ex}

\noindent [PSU20b]\, D. Popovici, J. Stelzig, L. Ugarte --- {\it Higher-Page Bott-Chern and Aeppli Cohomologies and Applications} --- J. reine angew. Math., doi: 10.1515/crelle-2021-0014.

\vspace{1ex}

\noindent [PSU20c]\, D. Popovici, J. Stelzig, L. Ugarte ---  {\it Deformations of Higher-Page Analogues of $\partial\bar\partial$-Manifolds} --- Math. Z. (2021), https://doi.org/10.1007/s00209-021-02861-0.

\vspace{1ex}

\noindent [PU23]\,  D. Popovici, L. Ugarte --- {\it A Moment Map for the Space of Maps to a Balanced Manifold} --- arXiv e-print DG 2311.00485v1; Bull. Sci. Math. {\bf 212} (2026), 103869; 
doi : 10.1016/j.bulsci.2026.103869

\vspace{1ex}

\noindent [Rol09]\,S. Rollenske --- {\it Geometry of nilmanifolds with left-invariant complex structure and deformations in the large} --- Proc. Lond. Math. Soc. {\bf 99} (2009), 425--460.

\vspace{1ex}

\noindent [Tia87]\, G. Tian --- {\it Smoothness of the Universal Deformation Space of Compact Calabi-Yau Manifolds and Its Petersson-Weil Metric} --- Mathematical Aspects of String Theory (San Diego, 1986), Adv. Ser. Math. Phys. 1, World Sci. Publishing, Singapore (1987), 629--646.

\vspace{1ex}

\noindent [Tod89]\, A. N. Todorov --- {\it The Weil-Petersson Geometry of the Moduli Space of $SU(n\geq 3)$ (Calabi-Yau) Manifolds I} --- Comm. Math. Phys. {\bf 126} (1989), 325-346.

\vspace{1ex}

\noindent [Tol25]\, 
A. Tolcachier --- {\it Six-dimensional complex solvmanifolds with non-invariant trivializing sections of their canonical bundle} --- arXiv e-print DG 2412.02325v2; 
Math. Nachr. {\bf 298} (2025), no. 8, 2626--2651, 
https://doi.org/10.1002/mana.70008



\vspace{3ex}

\noindent Graduate School of Mathematics    \hfill Institut de Math\'ematiques de Toulouse

\noindent Nagoya University, Furocho, Chikusaku,    \hfill  Universit\'e de Toulouse

\noindent Nagoya, Japan, 464-8602  \hfill  118 route de Narbonne, 31062 Toulouse, France

\noindent  Email: kasuya@math.nagoya-u.ac.jp                \hfill     Email: popovici@math.univ-toulouse.fr

\vspace{2ex}

\noindent and

\vspace{2ex}

\noindent Departamento de Matem\'aticas\,-\,I.U.M.A., Universidad de Zaragoza,

\noindent Campus Plaza San Francisco, 50009 Zaragoza, Spain

\noindent Email: ugarte@unizar.es

\end{document}